\documentclass[aap]{imsart}
\RequirePackage{amsmath,amsthm,amsfonts,amssymb,mathtools,mathrsfs}
\RequirePackage[numbers]{natbib}
\usepackage{microtype}
\usepackage{booktabs,array,graphicx}
\usepackage{xcolor}
\usepackage{xurl}
\usepackage[colorlinks=true,linkcolor=blue!45!black,citecolor=blue!45!black,urlcolor=blue!45!black]{hyperref}
\startlocaldefs
\theoremstyle{plain}
\newtheorem{theorem}{Theorem}[section]
\newtheorem{lemma}[theorem]{Lemma}
\newtheorem{proposition}[theorem]{Proposition}
\newtheorem{corollary}[theorem]{Corollary}
\theoremstyle{definition}

\newcommand{\E}{\mathbb E}
\newcommand{\R}{\mathbb R}
\newcommand{\Prb}{\mathbb P}

\numberwithin{equation}{section}
\allowdisplaybreaks[2]
\endlocaldefs
\begin{document}
\begin{frontmatter}
\title{Large deviations and conditional laws for a time-integrated logarithm of averaged Gaussian squares}
\runtitle{Large deviations for integrated Gaussian powers}
\runauthor{H. Xu and H. Zhao}
\begin{aug}
\author[A]{\fnms{Haoran}~\snm{Xu}\ead[label=e1]{20244007005@stu.suda.edu.cn}}
\author[B]{\fnms{Huazhong}~\snm{Zhao}\ead[label=e2]{1453785662@qq.com}}
\address[A]{School of Mathematical Sciences, Soochow University\printead[presep={,\ }]{e1}}
\address[B]{China Telecom Corporation Limited, Suzhou Branch\printead[presep={,\ }]{e2}}
\end{aug}
\hypersetup{pdfauthor={Haoran Xu and Huazhong Zhao},pdftitle={Large deviations and conditional laws for a time-integrated logarithm of averaged Gaussian squares}}
\begin{abstract}
Systems with identical instantaneous power distributions can have very different probabilities of insufficient service over a fixed interval. We study the lower tail of the time-integrated logarithm of a weighted mean of $d$ independent squared Gaussian paths as $d\to\infty$. A global entropy inequality resolves the nonconvex service constraint, identifies the unique Gaussian minimizer and controls the nonlinear likelihood weight far from that minimizer. Combined with a joint local limit, it yields sharp probabilities at fixed internal thresholds and conditional path laws: a Fredholm determinant enters the leading factor, and tangential Gaussian fluctuations are independent of the exponential service deficit. The approximations extend through moderate deviations to the central scale. For stationary Ornstein--Uhlenbeck paths, cancellation in the radial likelihood expansion yields relative error $O(d^{-1})$ on compact positive parameter families with proportions and saddles separated from zero. Genuine temporal covariance attenuation by correlation kernels for a type of positive limiting proportion increases the lower-tail rate, yet the finite-sample outage distributions cross near typical service. Thus one-time power distributions alone do not determine reliability.
\end{abstract}
\begin{keyword}[class=MSC]
\kwdgroup[type=primary]{\kwd{60F10}\kwd{60G15}}
\kwdgroup[type=secondary]{\kwd{60F05}}
\end{keyword}
\begin{keyword}
\kwd{Gaussian processes}\kwd{large deviations}\kwd{conditional limit theorem}\kwd{sharp asymptotics}\kwd{covariance comparison}
\end{keyword}
\end{frontmatter}

\section{Introduction}
\label{sec:introduction}

A high mean service rate does not ensure delivery by a deadline: in fading channels, a sustained period of low received power can exhaust the available time for transmission \cite{WuNegi2003}. With many independent modes, instantaneous power concentrates, but the persistence of fluctuations still matters. How small is the shortfall probability, which trajectories realize the shortfall, and how does temporal dependence affect it?

For $1\le k\le m$, let $X_{ki}$, $1\le i\le n_k$, be independent copies of a centered continuous real Gaussian process with law $P_k$ and covariance $C_{k,0}$ on $[0,T]$. Different types are independent, $d=\sum_k n_k$, and $p_k>0$ is the power weight. We study
\begin{equation}
 \begin{aligned}
 q_d(t)&=\frac1d\sum_{k=1}^m\sum_{i=1}^{n_k}p_kX_{ki}(t)^2,\\
 S_d&=\int_0^T\log(1+\rho q_d(t))\,dt,
 \qquad P_d(b)=\Prb(S_d\le b),
 \end{aligned}
 \label{eq:original-event}
\end{equation}
with $\rho,T>0$ fixed and $d\to\infty$. Covariance matrices at distinct times are assumed positive definite; stationarity and the Markov property are not required. Write $\pi_{k,d}=n_k/d$ and, in Section~\ref{sec:continuous-exponent}, assume $\pi_{k,d}\to\pi_k>0$. The typical service is
\[
 b_0=\int_0^T\log\left(1+\rho\sum_k\pi_kp_k C_{k,0}(t,t)\right)dt.
\]
Our interest is $0<b<b_0$. Models with identical covariance diagonals have the same one-time power distributions and mean service at every $d$, but can have very different lower tails.

For unit-variance paths and $n_k=2r_k$, pairing real paths gives circular complex gains $h_{kj}=(X_{k,2j-1}+\mathrm iX_{k,2j})/\sqrt2$ and $q_d=r^{-1}\sum_{k,j}p_k|h_{kj}|^2$, with $r=d/2$. The logarithm in \eqref{eq:original-event} is the instantaneous mutual information for a single receiver with channel knowledge and independent Gaussian inputs \cite{TseViswanath}; $\rho$ is the nominal total signal-to-noise ratio when $\sum_k\pi_{k,d}p_k=1$. Shi et al.~\cite{ShiMaYangTamXia2017} compute outage for finitely many correlated Nakagami fading rounds and derive high-signal-to-noise asymptotics; Kazakopoulos et al.~\cite{Kazakopoulos2011} derive multiantenna outage rates from random-matrix eigenvalue distributions. Here we increase the number of independent continuous paths at fixed signal-to-noise ratio and horizon.

\medskip\noindent\textit{Related work.}
Many-sources large deviations are familiar in queueing: Wischik \cite{Wischik2001} treats averages of discrete-time inputs; D\k{e}bicki and Mandjes \cite{DebickiMandjes2003} derive exact Gaussian overflow asymptotics, and Mandjes and van Uitert \cite{MandjesVanUitert2005} and Zubeldia and Mandjes \cite{ZubeldiaMandjes2021} use the reproducing-kernel rate function in networks. Unlike the Gaussian superposition in these works, squaring before averaging in \eqref{eq:original-event} makes the relevant change one of second moments, subject to a nonlinear constraint.

The closest precursor for our Ornstein--Uhlenbeck (OU) model is Donati-Martin, Rouault, Yor and Zani \cite{DonatiMartinRouaultYorZani2004}, who prove fixed-horizon large deviations for squared radial OU processes as the dimension grows. For $n$ independent unit-variance OU paths of rate $\kappa$, their mean square satisfies $dY=2\kappa(1-Y)dt+\sqrt{8\kappa Y/n}\,dW$. In addition to the dynamic action in their result, the stationary initial law $nY(0)\sim\chi_n^2$ incurs the cost $(c(0)-1-\log c(0))/2$ for a limiting profile $c$. Section~\ref{sec:ou-quantitative} identifies the minimum of this stationary action under the service constraint. The radial principle alone does not resolve that nonconvex optimization, identify the conditioned individual paths, or determine the probability prefactor; our approach also applies without a radial Markov representation.

Large deviations for quadratic Gaussian functionals over growing horizons are studied in \cite{BrycDembo1997,BercuGamboaLavielle2000,BercuRouault2002,ZamparoSemeraro2023}. For Gaussian fields, Adler, Blanchet and Liu \cite{AdlerBlanchetLiu2012} and Li and Liu \cite{LiLiu2015} develop efficient high-excursion simulation; Liu \cite{Liu2012} derives precise tails for exponential integrals, Liu and Xu \cite{LiuXu2014} obtain conditional fields in total variation, and Li, Liu and Xu \cite{LiLiuXu2016} study aggregate lognormal fields in small noise. In our regime the shortfall arises from a collective change of many path laws over a fixed interval.

Gibbs conditioning relates conditional laws to information projections \cite{Csiszar1984,CsiszarMatus2003}; see also conditioned subsequences \cite{BroniatowskiCaron2014}, localization under unbounded constraints \cite{Nam2020}, and rare-event mechanisms in fractional Brownian queues and moving averages \cite{AwadGlynn2009,ChakrabartySamorodnitsky2024}. Path-space large deviations are studied for interacting diffusions \cite{BaldassoPereiraReis2022} and in connection with stochastic control \cite{BackhoffLackerTangpi2020}. Chaintron, Conforti and Reygner \cite{ChaintronConfortiReygner2026} give conditional large deviations and necessary Gibbs optimality conditions for possibly nonconvex constraints, while Chaintron and Conforti \cite{ChaintronConforti2025} establish stability for convex marginal constraints. Applying the strengthened measure topology with squared-path control, as in the second-Wasserstein Sanov theorem \cite{WangWangWu2010}, would require all positive exponential moments of $\|X\|_\infty^2$, which Gaussian paths lack. Nonlinear deviation estimates of Chatterjee and Dembo \cite{ChatterjeeDembo2016} and Yan \cite{Yan2020} concern Bernoulli and compactly supported Banach-space variables. We overcome the nonconvexity by a global entropy inequality and use finite-time observations to establish the continuous exponent and conditioning.

For precise probabilities, scalar local limits \cite{BahadurRao1960,ChagantySethuraman1993}, curved boundaries in finite dimensions \cite{AndrianiBaldi,BarbeBroniatowski}, Markov random walks \cite{ChanLai2003} and projection geometry \cite{LiaoRamanan2024} explain how curvature affects prefactors. In Banach spaces, dominating-point results \cite{Kuelbs2000,EinmahlKuelbs2004} concern open convex sets; Laplace expansions also treat integral statistics \cite{Fatalov1999,Fatalov2004}. The moment assumptions of general expansions need care: Piterbarg and Fatalov \cite[Section~4]{PiterbargFatalov1995} assume all exponential norm moments of iid summands, and Albeverio and Liang \cite{AlbeverioLiang} require a square-exponential moment. For a unit vector $f$ with $\langle X,f\rangle$ nondegenerate, $\|X^2\|_2^2\ge T^{-1}|\langle X,f\rangle|^4$, which has no positive exponential moment. We use the global inequality to control distant fluctuations instead.

Small-noise stochastic-equation asymptotics also use fluctuations around optimal paths \cite{GrafkeSchaeferVandenEijnden2024} and Fredholm determinants \cite{SchorleppTongGrafkeStadler2023}. Schorlepp and Grafke \cite[Theorem~2.1]{SchorleppGrafke2026} treat smooth multiplicative noise with fixed initial data and a unique nondegenerate optimal noise. The stationary radial OU model has a square-root coefficient and random chi-square initial data. Our global argument includes non-Markov kernels, while an exact likelihood expansion yields a quantitative result for stationary OU paths.

\medskip\noindent\textit{Main results.}
We first minimize relative entropy under the nonlinear service constraint:
\[
 \inf\left\{\sum_k\pi_k D(Q_k\Vert P_k):
 \int_0^T\log\left(1+\rho\sum_k\pi_kp_k\E_{Q_k}X(t)^2\right)dt\le b\right\}.
\]
The feasible set is generally nonconvex. The logarithmic Jensen inequality underlying geodesic convexity of positive-definite matrices \cite{Wiesel2012,SraHosseini2015} gives a global service bound at a Gaussian tilt. Combined with the entropy difference and an entropy decomposition, it yields a quantitative relative-entropy gap for every competing finite-entropy path law, independent of time discretization. In particular, global optimality, uniqueness and curvature control are consequences of the model rather than assumptions on a candidate minimizer. The optimizer is given by
\begin{gather*}
 \frac{dQ_{k,g}}{dP_k}(X)=Z_k(g)^{-1}\exp\left\{-p_k\int_0^Tg(t)X(t)^2\,dt\right\},\\
 g(t)=\frac{\theta\rho}{1+\rho\sum_k\pi_kp_k\E_{Q_{k,g}}X(t)^2}.
\end{gather*}
A contraction in $\log g$ constructs the common potential, and a scalar $\theta>0$ selects the threshold. Theorems~\ref{thm:entropy-stability}, \ref{thm:continuous-exponent} and \ref{thm:conditioning} identify the entropy minimum with the outage exponent, prove conditional concentration of the power path in the uniform norm, and establish total variation convergence of any fixed collection of paths to independent Gaussian tilts.

Theorem~\ref{thm:sharp} gives the sharp probability at every fixed positive internal threshold:
\[
 P_d(b)\sim\frac{\det_F(I-K)^{-1/2}}{\theta\sqrt{2\pi dV}}e^{-dI(b)},
 \qquad V=\frac1\theta\langle1,K(I-K)^{-1}1\rangle .
\]
Here $I$ is the entropy minimum and $K$ is $\theta$ times the covariance of the linearized service-rate profile. The same entropy gap used to identify the optimum controls distant fluctuations under the nonlinear likelihood weight, despite the limited exponential moments of squared Gaussian paths. The curvature gap $\|K\|<1$ controls nearby fluctuations. Lemmas~\ref{lem:sqrt-residual} and \ref{lem:joint-local} combine these estimates with a joint local limit, giving the prefactor and a conditional tangential Gaussian independent of the exponential deficit $d(b-S_d)$, including convergence of conditional second moments. For stationary OU paths, the radial likelihood and boundary Jacobian yield a uniform $O(d^{-1})$ relative error by cancellation of the first correction (Theorem~\ref{thm:ou-relative}).

The application is a comparison of models with identical one-time distributions (Theorem~\ref{thm:attenuation}). Pointwise covariance attenuation by correlation kernels increases the lower-tail rate and decreases the Laplace transform and variance of the service, nonstrictly in general and strictly for genuine attenuation of a positively represented type. Equal means and different variances force the finite-sample distributions to cross. Shorter OU correlation times therefore reduce rare shortfalls but cannot decrease every failure probability.

The main results have different ranges: Theorems~\ref{thm:continuous-exponent}--\ref{thm:conditioning} use positive limiting proportions; Theorem~\ref{thm:sharp} allows zero limits at a positive limiting saddle, using the actual proportions in relative formulas. Theorem~\ref{thm:central-transition} and Corollary~\ref{cor:moderate-conditioning} cover the approach to zero saddle. The OU relative order holds on compact families with both proportions and multipliers separated from zero; its uniform constants are existential, not explicit finite-count error bounds. Proposition~\ref{prop:edgeworth} and the refined crossing displacement use fixed rational proportions and compatible counts. Further consequences include a quantile correction, a relative second-moment formula under the sufficient condition of Corollary~\ref{cor:is-moment}, a first central correction in the Lugannani--Rice form \cite{LugannaniRice1980}, and power-dispersion comparison. Section~\ref{sec:illustration} gives grid and continuous Riccati evaluation, finite-count continuous probability bounds, and numerical illustrations.

\section{Gaussian entropy projection and stability}
\label{sec:global-structure}

The multiplier $\theta$ selects a service level and the common potential $g$ adjusts the variance profiles. We first establish global optimality among all finite-entropy laws.

\subsection{Setting and a dimension-independent contraction}

Use the Gaussian path laws $P_k$ and powers $p_k$ of Section~\ref{sec:introduction}, with proportions $\pi_k>0$, $\sum_k\pi_k=1$, and put $H=\sum_k\pi_k\delta_{(P_k,p_k)}$. Write
\[
 m_0(t)=\sum_k\pi_kp_k C_{k,0}(t,t),\qquad
 M_k=\max_t C_{k,0}(t,t),\quad M=\max_t m_0(t)>0.
\]
Proportions may be the actual counts $n_k/d$; omit zero-count types. Constants hold throughout the simplex, with $M$ replaced if necessary by $\max_k p_kM_k$. For a positive $g\in C([0,T])$, set
\[
 \begin{aligned}
 Z_k(g)&=E_{P_k}\exp\{-p_k\int_0^T g(t)X(t)^2\,dt\},\\
 \frac{dQ_{k,g}}{dP_k}&=Z_k(g)^{-1}\exp\{-p_k\int_0^TgX^2\}.
 \end{aligned}
\]
Quadratic Gaussian transforms can be computed by Volterra filtering equations \cite{KleptsynaLeBreton2001}. Here we use the covariance operator on $L^2([0,T])$, with $M_g$ denoting multiplication by $g$. Diagonalizing the positive trace-class operator $p_k C_{k,0}^{1/2}M_g C_{k,0}^{1/2}$ gives independent scalar quadratic weights, each dividing its Gaussian variance by one plus twice its eigenvalue. The convergent products of normalizers (or dominated convergence through finite spectral projections) give the centered Gaussian tilt with
\[
 \begin{aligned}
 Z_k(g)&=\det_F(I+2p_k C_{k,0}^{1/2}M_g C_{k,0}^{1/2})^{-1/2},\\
 C_{k,g}&=C_{k,0}^{1/2}(I+2p_k C_{k,0}^{1/2}M_g C_{k,0}^{1/2})^{-1}C_{k,0}^{1/2}.
 \end{aligned}
\]
We use $(C_{k,0}^{-1}+2p_k M_g)^{-1}$ as shorthand for this bounded resolvent. Write
\[
 m_g(t)=\sum_k\pi_kp_k C_{k,g}(t,t),\quad w_g(t)=\frac{\rho}{1+\rho m_g(t)},\quad
 q_0=\frac{\rho M}{1+\rho M}<1.
\]

\begin{theorem}[Global logarithmic contraction]
\label{thm:global-contraction}
For every $\theta>0$, the map on $C([0,T])$
\[
 (\mathcal T_\theta h)(t)
 =\log(\theta\rho)-\log(1+\rho m_{e^h}(t))
\]
is a contraction in the supremum norm with factor at most $q_0$.
It has a unique fixed point $h_\theta=\log g_\theta$, and
\[
 \frac{\theta\rho}{1+\rho M}\le g_\theta(t)\le\theta\rho,
 \qquad
 g_\theta(t)(1+\rho m_{g_\theta}(t))=\theta\rho.
\]
For an arbitrary approximate potential $\widetilde h$,
\begin{equation}
 \|\widetilde h-h_\theta\|_\infty
 \le\frac{\|\widetilde h-\mathcal T_\theta\widetilde h\|_\infty}
 {1-q_0}.
 \label{eq:global-residual}
\end{equation}

\end{theorem}

\begin{proof}
The resolvent gives $0\preceq C_{k,g}\preceq C_{k,0}$, hence $0\le m_g\le m_0$. For $g+v>0$, the identity $C_{k,g+v}-C_{k,g}=-2p_k C_{k,g}M_vC_{k,g+v}$ bounds the kernel difference by $2p_kTM_k^2\|v\|_\infty$; substituting this bound back gives derivative $-2p_k C_{k,g}M_{\delta g}C_{k,g}$ with quadratic kernel-norm remainder. Thus
\[
 \delta m_g(t)=-2\sum_k\pi_kp_k^2\int_0^T C_{k,g}(t,s)^2\delta g(s)\,ds.
\]
Writing $C=C_{k,0}$ and $B=2p_kC^{1/2}M_gC^{1/2}$, the positive operator $C_{k,g}-2p_kC_{k,g}M_gC_{k,g}=C^{1/2}(I+B)^{-2}C^{1/2}$ has continuous kernel and nonnegative diagonal:
\begin{equation}
 2p_k\int_0^Tg(s)C_{k,g}(t,s)^2\,ds\le C_{k,g}(t,t).
 \label{eq:resolvent-row}
\end{equation}
Consequently $D\mathcal T_\theta(h)$ has nonnegative kernel
\[
 L_g(t,s)=2w_g(t)\sum_k\pi_kp_k^2 C_{k,g}(t,s)^2g(s)
\]
and row sums at most $w_g(t)m_g(t)\le q_0$. Integrating along a segment between two $h$'s proves the contraction. Banach's theorem gives the fixed point, the range of $\mathcal T_\theta$ gives its bounds, and the triangle inequality followed by absorption of the contraction term gives \eqref{eq:global-residual}.
\end{proof}

Starting at $h_0=\log(\theta\rho)$, each iteration evaluates the tilted covariances. With evaluation errors $\delta_j$ in supremum norm (zero for exact iteration),
\[
 \|h_N-h_\theta\|_\infty\le q_0^N\log(1+\rho M)
       +\sum_{j=0}^{N-1}q_0^{N-1-j}\delta_j.
\]

\begin{corollary}[Universal curvature gap]
\label{cor:global-gap}
At the fixed point, define the exact symmetric curvature kernel
\[
 K_\theta(t,s)
 =2\theta w_{g_\theta}(t)w_{g_\theta}(s)
       \sum_k\pi_kp_k^2C_{k,g_\theta}(t,s)^2.
\]
It is a positive trace-class operator and
\[
 0\preceq K_\theta,
 \qquad \|K_\theta\|_{\infty\to\infty}\le q_0,
 \qquad \|K_\theta\|_{2\to2}\le q_0<1.
\]
In particular the full Fredholm factor
$D=\det_F(I-K_\theta)^{-1/2}$ is well defined, with
\[
 \frac12\operatorname{tr}K_\theta
 \le\log D\le
 \frac{\operatorname{tr}K_\theta}{2(1-q_0)}.
\]
\end{corollary}

\begin{proof}
At the fixed point $g_\theta=\theta w_{g_\theta}$, so $K_\theta=L_{g_\theta}$ is symmetric. Squared covariance kernels are positive by the Schur product property; continuity and the diagonal integral give trace class. The row bound gives the $L^2$ bound by the Schur test. Apply $x\le-\log(1-x)\le x/(1-q_0)$ to the eigenvalues.
\end{proof}

\subsection{An entropy stability inequality for arbitrary path laws}

For competitors $Q=(Q_1,\ldots,Q_m)$ with $\mathcal H(Q):=\sum_k\pi_kD(Q_k\Vert P_k)<\infty$, put
\[
 m_Q(t)=\sum_k\pi_kp_kE_{Q_k}X(t)^2,\qquad
 F(Q)=\int_0^T\log(1+\rho m_Q(t))\,dt.
\]
The second moment includes any mean. Fernique's theorem \cite{Fernique1970} and the entropy inequality give $E_{Q_k}\|X\|_\infty^2<\infty$, hence continuity of $m_Q$. Let $Q^*=(Q_{k,g_\theta})_k$, $b_\theta=F(Q^*)$, and $\alpha_\theta=(1+2\theta\rho T\max_k p_kM_k)^{-1}>0$.

\begin{theorem}[Quantitative global information projection]
\label{thm:entropy-stability}
For every finite-entropy collection of path laws $Q$, including
non-Gaussian and noncentered laws,
\begin{equation}
 \mathcal H(Q)-\mathcal H(Q^*)
 +\theta\{F(Q)-b_\theta\}
 \ge
 \alpha_\theta\sum_k\pi_kD(Q_k\Vert Q_k^*).
 \label{eq:strong-projection}
\end{equation}
Consequently $Q^*$ uniquely minimizes $\mathcal H(Q)$ subject to
$F(Q)\le b_\theta$. Every such feasible law with
$\mathcal H(Q)\le\mathcal H(Q^*)+\varepsilon$ satisfies
\[
 \sum_k\pi_kD(Q_k\Vert Q_k^*)\le\varepsilon/\alpha_\theta,
 \qquad b_\theta-F(Q)\le\varepsilon/\theta.
\]
With total variation defined by $\sup_A|Q(A)-P(A)|$, one also has
\[
 \sum_k\pi_k\|Q_k-Q_k^*\|_{\rm TV}^2
 \le\varepsilon/(2\alpha_\theta).
\]
\end{theorem}

The coefficient in \eqref{eq:strong-projection} accounts for the curved constraint. The exact entropy identity is
\[
 \mathcal H(Q)-\mathcal H(Q^*)=\sum_k\pi_kD(Q_k\Vert Q_k^*)
       -\int_0^Tg_\theta(m_Q-m_{g_\theta}).
\]
With $u=(1+\rho m_Q)/(1+\rho m_{g_\theta})$, a boundary competitor satisfies $\int\log u=0$. Thus $\int g_\theta(m_Q-m_{g_\theta})=\theta\int(u-1-\log u)>0$ if the profiles differ.

\begin{proof}
Let distinct nodes $t_j$ have positive weights $a_j$ summing to $T$, and replace time integration by quadrature. Denote the prior and tilted covariance matrices by $C_k,S_k^*$, with $G=\operatorname{diag}(a_jg_j)$. The discrete contraction gives the unique solution
\[
 (S_k^*)^{-1}=C_k^{-1}+2p_kG,\qquad
 g_j=\frac{\theta\rho}{1+\rho\sum_k\pi_kp_k(S_k^*)_{jj}}.
\]
\emph{Gaussian competitors.} For centered Gaussian laws with covariances $S_k>0$, set
\[
 \begin{aligned}
 R_k&=(S_k^*)^{-1/2}S_k(S_k^*)^{-1/2},\quad
 A_k=(S_k^*)^{1/2}G(S_k^*)^{1/2},\\
 B_k&=(S_k^*)^{1/2}C_k^{-1}(S_k^*)^{1/2}=I-2p_kA_k.
 \end{aligned}
\]
At node $j$, write $v_{kj}=(S_k^*)^{1/2}e_j$, $m_j^*=\sum_k\pi_kp_k\|v_{kj}\|^2$. Weight the spectral probability measures of $\log R_k$ at $v_{kj}/\|v_{kj}\|$ by $\rho\pi_kp_k\|v_{kj}\|^2/(1+\rho m_j^*)$, and add mass $1/(1+\rho m_j^*)$ at zero. Jensen's inequality $\log\int e^x\,d\mu\ge\int x\,d\mu$ for the resulting probability measure gives
\[
 \log\frac{1+\rho\sum_k\pi_kp_kv_{kj}^{\mathsf T}R_kv_{kj}}{1+\rho m_j^*}
 \ge\frac{\rho\sum_k\pi_kp_kv_{kj}^{\mathsf T}(\log R_k)v_{kj}}{1+\rho m_j^*}.
\]
Multiplying by $\theta a_j$ and summing gives $\theta(F(S)-F(S^*))\ge\sum_k\pi_kp_k\operatorname{tr}(A_k\log R_k)$. The Gaussian entropy difference is
\[
 \mathcal H(S)-\mathcal H(S^*)=\tfrac12\sum_k\pi_k
 \{\operatorname{tr}(B_k(R_k-I))-\operatorname{tr}\log R_k\}.
\]
Since $B_k=I-2p_kA_k$, addition gives
\begin{equation}
 \begin{gathered}
 \mathcal H(S)-\mathcal H(S^*)+\theta(F(S)-F(S^*))\\
 \ge\tfrac12\sum_k\pi_k\operatorname{tr}\{B_k(R_k-I-\log R_k)\}.
 \end{gathered}
 \label{eq:matrix-entropy-gap}
\end{equation}
Both factors in the trace are positive semidefinite. By the eigenvalue formula and $g_j\le\theta\rho$,
\[
 \lambda_{\min}(B_k)=\frac1{1+2p_k\lambda_{\max}(C_k^{1/2}GC_k^{1/2})}
 \ge\frac1{1+2p_k\theta\rho T M_k}\ge\alpha_\theta,
\]
where the maximum eigenvalue is bounded by $\operatorname{tr}(C_k^{1/2}GC_k^{1/2})=\sum_ja_jg_j(C_k)_{jj}$. The right side of \eqref{eq:matrix-entropy-gap} therefore dominates $\alpha_\theta\sum_k\pi_kD(N(0,S_k)\Vert N(0,S_k^*))$.

\emph{Arbitrary vector laws.} For any finite-entropy law $q_k$, set $S_k=E_{q_k}XX^{\mathsf T}$ and $\gamma_k=N(0,S_k)$. Finite entropy relative to a nondegenerate Gaussian implies finite second moment and $S_k>0$. Equal second moments give
\[
 D(q_k\Vert N(0,C))=D(q_k\Vert\gamma_k)+D(\gamma_k\Vert N(0,C))
\]
for both $C=C_k$ and $C=S_k^*$. The service is unchanged by the replacement. Adding the entropy remainder to \eqref{eq:matrix-entropy-gap} and using $\alpha_\theta\le1$ proves the grid inequality.

\emph{Passage to path laws.} Take positive Riemann quadratures of vanishing mesh. Under the full path tilt $Q_{k,N}^*$ by $p_k\sum_ja_jg_{Nj}X(t_j)^2$, the conditional path law given the nodes is the same as under $P_k$. The entropy chain rule thus lifts the grid inequality, adding the identical conditional entropy on both sides with coefficients $1$ and $\alpha_\theta\le1$.

Extend the nodal potentials by $g_N(t)=\theta\rho/(1+\rho\sum_k\pi_kp_k C_{k,N}(t,t))$, where $C_{k,N}$ is the covariance of the full discrete tilt. The fixed-point bounds and covariance domination give equicontinuity:
\[
 |C_{k,N}(t,t)-C_{k,N}(s,s)|\le2\sqrt{M_k}
 \sqrt{E_{P_k}(X(t)-X(s))^2}.
\]
Along any uniformly convergent subsequence $g_N\to g>0$, discrete energies converge pathwise to $p_k\int gX^2$. Dominated convergence gives convergence of normalizers and covariance entries, so the limit solves the continuous equation. Uniqueness gives convergence of the whole sequence to $g_\theta$. The same compactness argument applies to convergent multipliers in a compact positive interval and convergent proportions.

For small $a>0$, $aE_{Q_k}\|X\|_\infty^2\le D(Q_k\Vert P_k)+\log E_{P_k}e^{a\|X\|_\infty^2}<\infty$. Thus the service quadrature and the energies in
\[
 D(Q_k\Vert Q_{k,N}^*)=D(Q_k\Vert P_k)+\log Z_{k,N}
       +p_kE_{Q_k}\sum_ja_jg_{Nj}X(t_j)^2
\]
converge, as does the entropy of $Q_{k,N}^*$ by covariance convergence. This yields \eqref{eq:strong-projection}. Rearrangement for feasible laws bounds the entropy excess below by the relative entropies plus $\theta(b_\theta-F(Q))$. Equality forces every relative entropy to vanish; Pinsker's inequality gives the total variation estimate.
\end{proof}

\subsection{The full saddle branch and its response}

Define $B(\theta)=b_\theta$ and $\Lambda(\theta)=\mathcal H(Q^*)+\theta B(\theta)$. The projection theorem and the scalar identity
$\theta\log(1+\rho x)=\inf_{a>0}\{ax+a/\rho-\theta-\theta\log(a/(\theta\rho))\}$
give the variational formulas
\[
 \begin{aligned}
 \Lambda(\theta)&=\inf_Q\{\mathcal H(Q)+\theta F(Q)\}\\
 &=\min_{g>0}\left\{\int_0^T\left[\frac g\rho-\theta
 -\theta\log\frac g{\theta\rho}\right]dt-\sum_k\pi_k\log Z_k(g)\right\}.
 \end{aligned}
\]
Indeed, minimizing jointly over $Q,g$ gives the second expression by the entropy variational formula, with $g_\theta$ attaining the minimum.

\begin{proposition}[Strict service response]
\label{prop:global-response}
For $\theta>0$,
\[
 \begin{gathered}
 u_\theta:=\frac{\partial\log g_\theta}{\partial\log\theta}
 =(I-K_\theta)^{-1}1,
 \qquad 1\le u_\theta\le\frac1{1-q_0},\\
 \Lambda'(\theta)=B(\theta),\qquad B'(\theta)=-V(\theta),\\
 V(\theta)=\frac1\theta
 \langle1,K_\theta(I-K_\theta)^{-1}1\rangle>0.
 \end{gathered}
\]
The map $B$ decreases continuously from
$\int_0^T\log(1+\rho m_0(t))dt$ to zero as $\theta$ increases
from zero to infinity. Every strictly internal threshold therefore has
one saddle. Its variational cost
\[
 I(b)=\inf_{F(Q)\le b}\mathcal H(Q)
 =\Lambda(\theta_b)-\theta_bb
\]
satisfies $I'(b)=-\theta_b$ and $I''(b)=1/V(\theta_b)>0$. It increases from zero to infinity as $b$ decreases from the typical service to zero.
\end{proposition}

\begin{proof}
The fixed-point derivative in $h$ is $K_\theta$; the contraction permits implicit differentiation in $\log\theta$, giving $(I-K_\theta)u_\theta=1$. The nonnegative Neumann series gives the bounds. Since $\partial_\theta g_\theta=w_{g_\theta}u_\theta$, the covariance derivative gives $B'(\theta)=\int w_{g_\theta}\partial_\theta m_{g_\theta}=-\theta^{-1}\langle1,K_\theta u_\theta\rangle$. The kernel is nonnegative and nonzero, so $V>0$. Differentiating the minimized objective by the envelope inequalities gives $\Lambda'=B$.

As $\theta\downarrow0$, $g_\theta\to0$ uniformly and the prior covariance is recovered. As $\theta\to\infty$, $\min g_\theta\to\infty$. For the covariance eigenvalues $\lambda_{kj}$, resolvent ordering gives
\[
 \int_0^T C_{k,g_\theta}(t,t)dt\le
 \sum_j\frac{\lambda_{kj}}{1+2p_k(\min g_\theta)\lambda_{kj}}\longrightarrow0
\]
by summable domination. Since $\log(1+x)\le x$, $B(\theta)\to0$. For a positive-weight type choose a unit $f\in L^2$ with $v_0=\langle f,C_{k,0}f\rangle>0$. Then $v_\theta=\langle f,C_{k,g_\theta}f\rangle\to0$ and data processing under $X\mapsto\langle f,X\rangle$ gives
\[
 D(Q_{k,g_\theta}\Vert P_k)\ge\tfrac12(v_\theta/v_0-1-\log(v_\theta/v_0))\longrightarrow\infty .
\]
The projection theorem gives the formula for $I$, the inverse-function theorem gives its derivatives, and the variational formula gives $I(b)=\sup_{\theta\ge0}(\Lambda(\theta)-\theta b)$: each feasible $Q$ bounds this expression above by $\mathcal H(Q)$, and the saddle attains equality. At $\theta=0$ the prior is optimal, with zero entropy.
\end{proof}

\section{Identification of the continuous reliability exponent}
\label{sec:continuous-exponent}

We identify the variational cost with the continuous-event exponent using finite-dimensional Gaussian large deviations and an exponentially good quadrature approximation.

\begin{lemma}[Gaussian modulus exponential moments]
\label{lem:modulus-mgf}
Let $X$ be a centered Gaussian random element of $C([0,T])$, and set
\[
 A=\|X\|_\infty,\qquad
 B_\delta=\sup_{|t-s|\le\delta}|X(t)-X(s)|.
\]
For every fixed $\lambda>0$,
\begin{equation}
 E\exp(\lambda A B_\delta)\longrightarrow1
 \quad\text{as }\delta\downarrow0.
 \label{eq:modulus-mgf}
\end{equation}
\end{lemma}

\begin{proof}
Continuity gives $B_\delta\to0$ almost surely and $B_\delta\le2A$. By Fernique \cite{Fernique1970}, $E e^{cA^2}<\infty$ for some $c>0$, hence $\mu_\delta=EB_\delta\to0$. Also $\sigma_\delta^2=\sup_{|t-s|\le\delta}E(X(t)-X(s))^2\to0$. Borell--TIS \cite[Theorem~2.1.1]{AdlerTaylor} (from Gaussian isoperimetry \cite{Borell1975}) applied to the increments and their negatives gives $\Pr(B_\delta>\mu_\delta+u)\le e^{-u^2/(2\sigma_\delta^2)}$ for $u>0$; the case $\sigma_\delta=0$ is immediate. With $Y=(B_\delta-\mu_\delta)_+$, tail integration gives $E e^{2aY^2}\le1+\int_0^\infty4au e^{2au^2-u^2/(2\sigma_\delta^2)}du=(1-4a\sigma_\delta^2)^{-1}$ for $4a\sigma_\delta^2<1$. Since $B_\delta^2\le2\mu_\delta^2+2Y^2$,
\begin{equation}
 E e^{aB_\delta^2}\le\frac{e^{2a\mu_\delta^2}}{1-4a\sigma_\delta^2}\longrightarrow1.
 \label{eq:modulus-square-mgf}
\end{equation}
Choose $a_0>0$ with $2a_0<c$. Young's and Cauchy--Schwarz inequalities give
\[
 E e^{2\lambda A B_\delta}
 \le (E e^{2a_0 A^2})^{1/2}(E e^{(2\lambda^2/a_0)B_\delta^2})^{1/2}.
\]
By \eqref{eq:modulus-square-mgf} the right side is bounded for small $\delta$. Uniform integrability and almost sure convergence prove \eqref{eq:modulus-mgf}.
\end{proof}

\begin{theorem}[Continuous Gaussian service exponent]
\label{thm:continuous-exponent}
Fix the type laws, powers, $T,\rho$ of Section~\ref{sec:global-structure}, with integer counts satisfying $\sum_kn_k(d)=d$ and $\pi_{k,d}=n_k(d)/d\to\pi_k>0$. Use the service \eqref{eq:original-event} and typical service $b_0$ of Section~\ref{sec:introduction}. For $0<b<b_0$,
\begin{equation}
 \lim_{d\to\infty}-\frac1d\log\Pr(S_d\le b)=I_H(b),
 \label{eq:true-exponent}
\end{equation}
where $I_H$ is the exact information-projection cost in
Proposition~\ref{prop:global-response}. Moreover, with the actual
$H_d=\sum_k\pi_{k,d}\delta_{(P_k,p_k)}$,
\[
 -\frac1d\log\Pr(S_d\le b)-I_{H_d}(b)\longrightarrow0.
\]

\end{theorem}

\begin{proof}
Partition the interval into cells of lengths $a_{Nj}$, representatives $t_{Nj}$ and mesh $\delta_N\to0$. Define $S_{d,N}=\sum_ja_{Nj}\log(1+\rho q_d(t_{Nj}))$. At fixed $N$, the sampled vectors $X_{ki,N}$ have positive-definite covariance matrices $C_{k,N}$. The sample second moments $U_{k,d}=n_k(d)^{-1}\sum_iX_{ki,N}X_{ki,N}^{\mathsf T}$ obey Cram\'er's theorem \cite[Corollary~6.1.6]{DemboZeitouni2010} since the quadratic moment generating function is finite near zero. The rate, for $U>0$, is
\[
 J_{k,N}(U)=\tfrac12\{\operatorname{tr}(C_{k,N}^{-1}U)-\log\det(C_{k,N}^{-1}U)-N\},
\]
and is infinite otherwise. Indeed the Legendre objective $\operatorname{tr}(LU)+\tfrac12\log\det(I-2C_{k,N}^{1/2}LC_{k,N}^{1/2})$ is maximized at $L=(C_{k,N}^{-1}-U^{-1})/2$ in the positive-definite case. Otherwise take a unit $v$ with $v^{\mathsf T}Uv\le0$ and $L=-svv^{\mathsf T}$; the objective diverges as $s\to\infty$.

Independence and $\pi_{k,d}\to\pi_k>0$ give the joint rate $\sum_k\pi_kJ_{k,N}(U_k)$ at speed $d$. The tuple is exponentially tight by a positive exponential trace moment. Replacing actual by limiting proportions changes the grid service by at most $C\max_k|\pi_{k,d}-\pi_k|\sum_k\operatorname{tr}U_k$, hence is superexponentially negligible by Chernoff's inequality. The contraction principle and the grid version of Theorem~\ref{thm:entropy-stability} identify the constrained minimum as $I_{H,N}(b)$ for large $N$. Its infima on $F_N<b$ and $F_N\le b$ agree by continuity of the saddle branch from below. Thus
\begin{equation}
 \lim_{d\to\infty}-d^{-1}\log\Pr(S_{d,N}\le b)=I_{H,N}(b).
 \label{eq:finite-grid-exponent}
\end{equation}
The discrete-tilt convergence in Theorem~\ref{thm:entropy-stability} gives $B_N(\theta)\to B(\theta)$ and convergence of entropies at fixed $\theta$. Bracket $\theta_b$ by $\theta_-<\theta_b<\theta_+$ with $B(\theta_-)>b>B(\theta_+)$. For large $N$ the grid root belongs to this bracket; shrinking it gives $\theta_{N,b}\to\theta_b$. The same compactness argument at these multipliers yields
\begin{equation}
 I_{H,N}(b)\longrightarrow I_H(b).
 \label{eq:grid-rate-limit}
\end{equation}
For the continuous event, writing $A_{ki},B_{ki,\delta}$ as in Lemma~\ref{lem:modulus-mgf},
\[
 |S_d-S_{d,N}|\le\frac{2\rho T}{d}\sum_{k,i}p_k A_{ki}B_{ki,\delta_N}.
\]
Hence for any $\varepsilon,\lambda>0$ the upper exponential rate of $\{|S_d-S_{d,N}|>\varepsilon\}$ is at most $-\lambda\varepsilon+\sum_k\pi_k\log E_{P_k}e^{2\lambda\rho Tp_kAB_{\delta_N}}$. Since the sum tends to zero and $\lambda$ is arbitrary,
\begin{equation}
 \lim_{N\to\infty}\limsup_{d\to\infty}\frac1d
       \log\Pr(|S_d-S_{d,N}|>\varepsilon)=-\infty.
 \label{eq:superexp-grid}
\end{equation}
Sandwich the continuous probability by grid probabilities at $b\pm\varepsilon$ plus or minus the error event. For large $N$ its decay exponent exceeds $I_{H,N}(b-\varepsilon)$, so subtraction preserves the lower bound's exponent. Equations \eqref{eq:finite-grid-exponent}--\eqref{eq:superexp-grid}, followed by $\varepsilon\downarrow0$, give \eqref{eq:true-exponent}. Finally contraction gives continuous dependence on the finite vector of proportions, and the same saddle bracketing gives $I_{H_d}(b)\to I_H(b)$.
\end{proof}

\subsection{Conditional power profiles and path laws}

We next identify the paths realizing the continuous rare event.

\begin{theorem}[Conditional concentration and asymptotic independence]
\label{thm:conditioning}
Under the assumptions of Theorem~\ref{thm:continuous-exponent}, fix $0<b<b_0$. Let $g=g_{\theta_b}$ and $m_g$ be formed with the limiting proportions, and write
\[
 \mathcal R_d=\mathcal L((X_{ki})_{k,i}\mid S_d\le b),\qquad
 \mathcal Q_d^*=\bigotimes_k Q_{k,g}^{\otimes n_k(d)}.
\]
For every $\eta>0$,
\begin{equation}
 \limsup_{d\to\infty}\frac1d\log\Prb\big(\|q_d-m_g\|_\infty\ge\eta\mid S_d\le b\big)<0.
 \label{eq:conditional-profile}
\end{equation}
Moreover,
\begin{equation}
 \frac1d D(\mathcal R_d\Vert\mathcal Q_d^*)\longrightarrow0.
 \label{eq:conditional-entropy}
\end{equation}
Consequently, for any fixed nonnegative integers $r_k$, the conditional law of the first $r_k$ paths of each type converges in total variation to $\bigotimes_k Q_{k,g}^{\otimes r_k}$. The empirical law within each type converges weakly in probability, under the conditioning, to $Q_{k,g}$.
\end{theorem}

\begin{proof}
\emph{Entropy cost of a power deviation.} Use the same grids. For $c$ in a compact neighborhood of $b$, the saddle $\theta_{N,c}\le\theta_+$ uniformly on fine grids by bracketing. Write $S_{k,N,c}^*$ for the optimal matrices, and set
\[
 a_*=(1+2\theta_+\rho T\max_k p_kM_k)^{-1},\qquad
 \mathcal D_N(U)=\sum_k\pi_k D(N(0,U_k)\Vert N(0,S_{k,N,c}^*)).
\]
For grid service at most $c$ (limiting proportions), the entropy gap gives
\begin{equation}
 \sum_k\pi_kJ_{k,N}(U_k)\ge I_{H,N}(c)+a_*\mathcal D_N(U).
 \label{eq:conditional-grid-gap}
\end{equation}
At node $j$ the centered variable $p_k(X_j^2-(S_{k,N,c}^*)_{jj})$ under the tilt has law $a(Z^2-1)$, $Z$ standard normal, $0\le a\le\max_kp_kM_k$. Its log transform $-ua-\frac12\log(1-2ua)$ is bounded by $Cu^2$ for $|u|\le u_0$, uniformly in all indices. The entropy inequality with both signs gives, for $0<u\le u_0$,
\begin{equation}
 \left|\sum_k\pi_kp_k((U_k)_{jj}-(S_{k,N,c}^*)_{jj})\right|
 \le\frac{\mathcal D_N(U)}{u}+Cu.
 \label{eq:conditioning-node-entropy}
\end{equation}
For $\eta>0$ choose $u$ with $Cu\le\eta/8$. A deviation of at least $\eta/4$ at a node costs $\mathcal D_N(U)\ge u\eta/8$, giving the joint rate lower bound $I_{H,N}(c)+\gamma$, $\gamma=a_*u\eta/8>0$ independent of mesh. Singular matrices have infinite rate.

\emph{Conditional concentration.} For the power modulus at mesh $\delta$, use the bound for $|S_d-S_{d,N}|$ above with $\rho T$ omitted. As in \eqref{eq:superexp-grid}, its exceeding any positive cutoff is superexponentially negligible in the successive limits. On a fixed grid, replacing actual by limiting proportions in the node powers and service also has superexponentially small error by the trace moment bound used above. Choose $c=b+3\varepsilon<b_0$ with $I_H(c)>I_H(b)-\gamma/4$ and $\|m_{g_{\theta_c}}-m_g\|_\infty<\eta/8$. Take a grid fine enough that its optimal profile at $c$ differs from $m_{g_{\theta_c}}$ by less than $\eta/8$, $|I_{H,N}(c)-I_H(c)|<\gamma/4$, and the cell modulus of $m_g$ is less than $\eta/8$. Require also upper exponential rates less than $-I_H(b)-\gamma$ for service error exceeding $\varepsilon$ and power modulus exceeding $\eta/8$. Outside these errors and proportion errors with the same cutoffs, $S_d\le b$ implies grid service at most $c$; a power deviation $\|q_d-m_g\|_\infty\ge\eta$ implies deviation at some node of at least $\eta/4$ after subtraction of the cell, proportion and profile errors. The closed-set LDP upper bound and \eqref{eq:conditional-grid-gap} give
\[
 \limsup_{d\to\infty}d^{-1}\log\Prb(S_d\le b,\ \|q_d-m_g\|_\infty\ge\eta)
 \le-I_H(b)-\gamma/2.
\]
Division by $P_d(b)$ using Theorem~\ref{thm:continuous-exponent} proves \eqref{eq:conditional-profile}.

For the entropy assertion, $\|q_d\|_\infty\le d^{-1}\sum_{k,i}p_k\|X_{ki}\|_\infty^2$ gives $\Prb(\|q_d\|_\infty>L)\le e^{-d(aL-C_0)}$ for some $a>0,C_0<\infty$ by Fernique. Conditioning changes the exponent by at most $d(I_H(b)+1)$ for large $d$. The integrated tail and \eqref{eq:conditional-profile} give $\E_{\mathcal R_d}\|q_d-m_g\|_\infty\to0$. The density of the conditioned law yields
\begin{align*}
 d^{-1}D(\mathcal R_d\Vert\mathcal Q_d^*)
 &= -d^{-1}\log P_d(b)+\sum_k\pi_{k,d}\log Z_k(g)
       +\E_{\mathcal R_d}\int_0^Tgq_d\\
 &\longrightarrow I_H(b)+\sum_k\pi_k\log Z_k(g)+\int_0^Tgm_g=0,
\end{align*}
by the tilt's entropy identity.

\emph{Tagged paths.} For nonzero $(r_k)$ take $\min_{k:r_k>0}\lfloor n_k/r_k\rfloor$ disjoint blocks of this composition, of order $d$. By the entropy chain rule and convexity, total entropy relative to the product dominates the sum of block marginal entropies. By exchangeability all block marginals coincide; their entropy relative to $\bigotimes_k Q_{k,g}^{\otimes r_k}$ tends to zero by \eqref{eq:conditional-entropy}. Pinsker gives total variation convergence. One- and two-path convergence give convergence of the mean and vanishing variance of the empirical average of each bounded continuous test function. A countable convergence-determining class on the Polish path space gives weak convergence in probability.
\end{proof}

\section{Sharp probabilities and conditional fluctuations}
\label{sec:sharp}

Precise deviations for sample covariances involve normal approximations \cite{FangKoike2022} and scalar extreme tails \cite{HeinyMikoschYslas2021}. Here the prefactor requires the mass of a narrow strip along a nonlinear boundary in path space. Throughout this section use actual proportions $\pi_{k,d}$, omitting types of zero count; replacing them by their limits need not preserve a relative approximation.

For a given $\theta>0$, abbreviate $b=B(\theta)$, $g=g_\theta$, $w=g/\theta$, and let $K=K_\theta$ be the kernel in Corollary~\ref{cor:global-gap}. Define
\begin{equation}
 \begin{gathered}
 D=\det_F(I-K)^{-1/2},\qquad
 V=\frac1\theta\langle1,K(I-K)^{-1}1\rangle,\\
 P_{0,d}(b)=\frac{D}{\theta\sqrt{2\pi dV}}e^{-dI(b)}.
 \end{gathered}
 \label{eq:effective-centre}
\end{equation}
Inner products and operator norms below are on $L^2([0,T])$ unless specified otherwise.

\begin{theorem}[Sharp asymptotics and conditional fluctuations]
\label{thm:sharp}
Fix the finite list of Gaussian types and powers, $T$ and $\rho$ as in Section~\ref{sec:global-structure}. Suppose $\pi_{k,d}\to\pi_k$ (zero limiting proportions are allowed) and $\theta_d\to\theta>0$. Use the actual saddle $g_d=g_{\theta_d}$ and threshold $b_d=B(\theta_d)$ at dimension $d$. Then
\begin{equation}
 P_d(b_d)=P_{0,d}(b_d)(1+o(1)).
 \label{eq:sharp-general}
\end{equation}
Let $K_\theta$ denote the curvature kernel at the limiting proportions and multiplier, and set
\[
 \Sigma=\theta^{-1}K_\theta(I-K_\theta)^{-1},\qquad V=\langle1,\Sigma1\rangle,
 \qquad \Sigma_\perp=\Sigma-\frac{(\Sigma1)\otimes(\Sigma1)}{V}.
\]
The notation $v\otimes v$ denotes $f\mapsto\langle v,f\rangle v$. Under $S_d\le b_d$,
\begin{equation}
 \left(\frac{g_d}{\theta_d}\sqrt d(q_d-m_{g_d}),\ d(b_d-S_d)\right)
 \ \Longrightarrow\ (Z_\perp,E_\theta)
 \label{eq:conditional-fluctuation}
\end{equation}
in $L^2([0,T])\times\mathbb R$. The random element $Z_\perp$ is centered Gaussian with covariance $\Sigma_\perp$, and $E_\theta$ is independent of $Z_\perp$ and exponential with rate $\theta$. For every fixed $f\in L^2([0,T])$, the conditional moments through order two of the pair
$(\langle f,w_d U_d\rangle,\theta_d d(b_d-S_d))$ converge to those of
$(\langle f,Z_\perp\rangle,\theta E_\theta)$, where $U_d=\sqrt d(q_d-m_{g_d})$ and $w_d=g_d/\theta_d$. In particular the limiting mean of the scaled deficit is one and its limiting covariance with each such projection is zero.
\end{theorem}

For fixed $0<b<b_0$, actual saddles converge to $\theta_b>0$ by Proposition~\ref{prop:global-response}. The proof combines ordinary moment and interval estimates, exponential control of the boundary remainder, and a joint local limit in a narrow normal strip (Lemmas~\ref{lem:central-half-estimates}--\ref{lem:joint-local}).

\subsection{The exact change of law and uniform moment bounds}

For notational convenience suppress the dependence of the saddle on $d$. Let $Q$ be the product of the tilted path laws and, for a coordinate of type $k$, put $Y_i(t)=p_k(X_i(t)^2-C_{k,g}(t,t))$. Write
\begin{align*}
 U_d(t)&=d^{-1/2}\sum_i Y_i(t)=\sqrt d(q_d(t)-m_g(t)),&
 L_d&=\int_0^T w(t)U_d(t)dt,\\
 y_d(t)&=w(t)U_d(t)/\sqrt d,&
 R_d&=d\int_0^T[y_d(t)-\log(1+y_d(t))]dt.
\end{align*}
The ratio $1+y_d=(1+\rho q_d)/(1+\rho m_g)$ is bounded below by $(1+\rho M)^{-1}$. Therefore $R_d\ge0$, and $d(S_d-b)=\sqrt d L_d-R_d$. The product likelihood is exactly
\begin{equation}
 \frac{d\Prb}{dQ}=\exp(-dI+\theta\sqrt d L_d),\qquad
 e^{dI}P_d(b)=\E_Q[e^{\theta\sqrt d L_d}\mathbf1_{\{\sqrt d L_d\le R_d\}}].
 \label{eq:sharp-change-of-law}
\end{equation}
The quadratic approximation to the boundary is
\begin{equation}
 R_{2,d}=\frac12\int_0^T(wU_d)^2dt,\qquad
 |R_d-R_{2,d}|\le \frac{C}{\sqrt d}\int_0^T|U_d(t)|^3dt.
 \label{eq:sharp-remainder}
\end{equation}
This follows by Taylor's theorem on $y\ge -1+(1+\rho M)^{-1}$; the third derivative of $\log(1+y)$ is bounded there. In the uniform estimates below $M$ may always be replaced by $\max_k p_kM_k$.

Split the coordinates into disjoint sets $A,B$ of sizes $\lfloor d/2\rfloor,\lceil d/2\rceil$, assigning the count of each type as evenly as possible between the halves (distribute odd remaining coordinates to balance the total sizes). Let $U_A,U_B,L_A,L_B$ be the corresponding partial sums with the same normalization $\sqrt d$. For $j=A,B$ define
\[
 \psi_d(v)=\min(v^2,\sqrt d\,|v|),\qquad J_j=\int_0^T\psi_d(U_j(t))dt.
\]

\begin{lemma}[Two-part estimates]
\label{lem:central-half-estimates}
For every fixed finite upper bound on $\theta$, there are positive constants $C_R,c,C$, depending only on this bound and the fixed data, such that uniformly over actual type counts, positive $\theta$ in this range and $d\ge4$,
\begin{equation}
 R_d\le C_R(J_A+J_B),\qquad \E_Q e^{cJ_j}\le C,
 \qquad \int_0^T\E_Q|U_j(t)|^3dt\le C.
 \label{eq:half-moments}
\end{equation}
For every interval $J\subset\mathbb R$,
\begin{equation}
 Q(L_j\in J)\le C(|J|+d^{-1/2}),\qquad
 \sup_z|Q(L_d/s\le z)-\Phi(z)|\le Cd^{-1/2},
 \label{eq:half-interval}
\end{equation}
where $\Phi$ is the standard normal CDF and $s^2=\operatorname{Var}_Q L_d=\theta^{-1}\langle1,K1\rangle$ stays bounded above and below by positive constants. The interval bound also holds for the partial sum over any block with at least $\delta d$ coordinates, for fixed $\delta>0$ and a constant depending on $\delta$.
\end{lemma}

\begin{proof}
On $y\ge-1+(1+\rho M)^{-1}$, $0\le y-\log(1+y)\le C\min(y^2,|y|)$ by Taylor expansion and linear growth. Thus $w\le\rho$ and $\psi_d(x+y)\le4(\psi_d(x)+\psi_d(y))$ give the first bound. At each time $Y_i$ has law $a_i(Z_i^2-1)$, $0\le a_i\le\max_k p_kM_k$, with independent standard normals. The quadratic transform gives $\log\E_Q e^{uU_j(t)}\le C_1u^2$ for $|u|\le c_1\sqrt d$. Chernoff \cite{BLM2013} with $|u|=c_2\min(v,\sqrt d)$ yields $Q(|U_j(t)|>v)\le2e^{-c_3\psi_d(v)}$ for $v\ge0$. Integrating gives bounded third moments and an exponential moment of $\psi_d(U_j(t))$ at a small positive argument; Jensen for $dt/T$ bounds $\E_Q e^{cJ_j}$.

The unscaled scalar summands $\int wY_i$ have bounded third moments and variance $2p_k^2\operatorname{tr}((C_{k,g}^{1/2}M_w C_{k,g}^{1/2})^2)$. The resolvent bounds $C_{k,g}$ between $C_{k,0}$ and a fixed positive multiple of it, and $w$ has a common positive lower bound. Hence the largest eigenvalue of $p_k C_{k,g}^{1/2}M_w C_{k,g}^{1/2}$ is uniformly bounded away from zero over types and saddles; the trace gives an upper variance bound. Berry--Esseen for independent summands \cite{ChenGoldsteinShao2011} gives normal error $O_\delta(d^{-1/2})$ for any block of at least $\delta d$ terms, with variance bounded above and below. Bounded normal densities give the interval estimates; Wick's identity gives $s^2=\theta^{-1}\langle1,K1\rangle$.
\end{proof}

\subsection{Uniform exponential moments of the boundary remainder}

The likelihood weight can amplify atypical profiles even under the tilt. A moment of order greater than one for $e^{\theta R_d}$ will suffice: we use the global entropy gap away from the saddle and a truncated Gaussian bound nearby, together with balanced blocks to retain interval estimates.

\begin{lemma}[Exponential moments of the nonlinear remainder]
\label{lem:sqrt-residual}
Consider a sequence of actual count vectors and exact saddles with convergent proportions (zero limits allowed) and $\theta_d\to\theta>0$. Take an upper bound $\theta_+$ for the multipliers, put
\[
 \alpha_*=(1+2\theta_+\rho T\max_k p_kM_k)^{-1}.
\]
For every $s>1$ with
\begin{equation}
 s(1-\alpha_*)<1,
 \label{eq:residual-moment-range}
\end{equation}
the full-sample residual $\mathcal T_d=\theta_dR_d$ satisfies
$\sup_d\E_Q e^{s\mathcal T_d}<\infty$.
In particular such an $s$ always exists. For sufficiently large $d$, the same bound holds for the residual of any one of a fixed number $r$ of nonempty blocks, with the count of each type distributed as evenly as possible among the blocks, evaluated at the full saddle. For a block $j$ of size $n_j$ let $q_j$ be its average weighted power, and put
\[
 \mathcal T_j=n_j\theta_d\int_0^T
 \phi\!\left(\frac{g_d(t)}{\rho\theta_d}(1+\rho q_j(t))\right)dt,\qquad
 \phi(z)=z-1-\log z.
\]
These independent block residuals satisfy
\begin{equation}
 \theta_d R_d\le\sum_{j=1}^r\mathcal T_j .
 \label{eq:residual-split}
\end{equation}
The bounds are uniform over compact positive multiplier intervals and all proportions if $s$ satisfies \eqref{eq:residual-moment-range} with a uniform strict gap. This uniformity also holds for stationary OU types on compact positive power and decay-rate families.
\end{lemma}

\begin{proof}
\emph{Quadratic transforms and the curvature bound.} First consider the full sample, suppressing $d$ in the saddle. Write
\[
 \mathcal R(q)=\theta\int_0^T\phi\!\left(\frac{g}{\rho\theta}(1+\rho q)\right)dt,\qquad
 \mathcal T_d=d\mathcal R(q_d).
\]
For $q\ge0$ the argument of $\phi$ is at least $a=(1+\rho\max_k p_kM_k)^{-1}$; $\phi$ is Lipschitz on $[a,\infty)$, and $\mathcal R(q)\le\int gq+C$. With $A_k=M_{\sqrt g}C_{k,0}M_{\sqrt g}$, the resolvent gives
\[
 M_{\sqrt g}C_{k,g}M_{\sqrt g}=A_k(I+2p_kA_k)^{-1},\quad
 2p_k\|M_{\sqrt g}C_{k,g}M_{\sqrt g}\|\le1-\alpha_*,
\]
since $\|A_k\|\le\theta_+\rho T M_k$. For $s_2>s$ satisfying \eqref{eq:residual-moment-range}, the quadratic transform gives
\begin{equation}
 \E_Q e^{s_2\mathcal T_d}\le e^{Cd}.
 \label{eq:crude-residual-moment}
\end{equation}
The eigenvalue sums are bounded and separated uniformly from the singularity. We also need an improved curvature bound. Writing $C=C_{k,0}$ and $B=2p_kC^{1/2}M_gC^{1/2}$, the spectral inequality $B(I+B)^{-1}\preceq(1-\alpha_*)I$ gives
\[
 2p_k C_{k,g}M_gC_{k,g}
 =C^{1/2}B(I+B)^{-2}C^{1/2}\preceq(1-\alpha_*)C_{k,g}.
\]
Test on normalized averages over shrinking intervals to pass this comparison to diagonals. The row argument of \eqref{eq:resolvent-row} gives $\|K_{\theta_d}\|\le(1-\alpha_*)q_0$ with a common $q_0<1$; thus $sq_*<1$ for $q_*=\limsup_d\|K_{\theta_d}\|$.

\emph{Away from the saddle.} We show that for each $\delta>0$ some $\gamma_\delta>0$ satisfies
\begin{equation}
 \limsup_d\frac1d\log\E_Q
 [e^{s\mathcal T_d}\mathbf1_{\{\|q_d-m_g\|_\infty\ge\delta\}}]\le-\gamma_\delta.
 \label{eq:weighted-far}
\end{equation}
First suppose the limiting proportions are positive. By Fernique and independence, the sup norm of $q_d$ has a tail with arbitrarily large exponential rate as its cutoff grows. For its modulus use $2d^{-1}\sum_i p_i\|X_i\|_\infty B_{i,h}$, as in Section~\ref{sec:continuous-exponent}. The tilt density is uniformly bounded by $Z_k(g)^{-1}$; for $G_h=\lambda AB_h$, $\E_{Q_{k,g}}(e^{G_h}-1)\le Z_k(g)^{-1}\E_{P_k}(e^{G_h}-1)\to0$ by Lemma~\ref{lem:modulus-mgf}. Thus exceeding any fixed modulus tolerance has arbitrarily negative rate on fine grids. H\"older with \eqref{eq:crude-residual-moment} gives the same conclusions with weight $e^{s\mathcal T_d}$.

On the complementary bounded-power, small-modulus set, use a fixed positive quadrature at its exact saddle $g_{N,d},m_{N,d}$ for the same multiplier and actual proportions, with product path tilt $Q_{N,d}$. The extensions and normalizers converge by Theorem~\ref{thm:entropy-stability}. The log density ratio is
\[
 \frac1d\log\frac{dQ}{dQ_{N,d}}
 =\sum_j a_jg_{N,d}(t_j)q_d(t_j)-\int_0^Tg_dq_d
       +\sum_k\pi_{k,d}\log\frac{Z_{k,N,d}}{Z_{k,g_d}}.
\]
Uniform bounds, node convergence and the modulus tolerance make this arbitrarily small, as well as the difference between $\mathcal R(q_d)$ and $\mathcal R_{N,d}(q)=\theta_d\sum_j a_j\phi(g_{N,d}(t_j)(1+\rho q_j)/(\rho\theta_d))$ at the sampled powers, using Lipschitz $\phi$. With sufficiently fine approximations the deviation $\|q_d-m_g\|_\infty\ge\delta$ implies a node deviation at least $\delta/2$.

At the limiting grid saddle, for positive matrices $U_k$ put
\[
 q_j=\sum_k\pi_kp_k(U_k)_{jj},\qquad
 \mathscr D_N(U)=\sum_k\pi_kD(N(0,U_k)\Vert N(0,S^*_{k,N})).
\]
The entropy identity and \eqref{eq:strong-projection} give $\mathscr D_N-\mathcal R_N\ge\alpha_*\mathscr D_N$. By \eqref{eq:conditioning-node-entropy}, a node deviation at least $\delta/3$ costs $\mathscr D_N\ge e_\delta>0$ independent of grid and positive proportions. Hence there
\[
 s\mathcal R_N(q)-\mathscr D_N(U)\le-[1-s(1-\alpha_*)]e_\delta<0.
\]
Under $Q_{N,d}$ the empirical node matrices have Cram\'er rate $\mathscr D_N$ at speed $d$: couple the varying covariances to their limits by applying the respective matrix square roots to common standard normals. The moment-matrix difference is bounded by a vanishing coefficient times empirical squared norms, so is superexponentially negligible. Positive limiting proportions and bounded node powers imply bounded traces of every $U_k$, by the bounds proportional to $1/(\pi_{k,d}p_k)$. On this compact set the residuals converge uniformly and the deviation implies a limiting deviation of at least $\delta/3$. Partition the bounded residual weight into level intervals of arbitrarily small width and apply the closed-set LDP upper bound. This gives the strict negative weighted rate above. Choose the grid likelihood and residual comparison errors together smaller than half the gap. This proves \eqref{eq:weighted-far}, with negative rate independent of any lower bound on the positive proportions.

For zero limiting proportions, discard $\|q_d\|_\infty>R$ by the weighted sup-norm tail. Add $\lceil\varepsilon d\rceil$ independent paths of every type at the original potential $g_d$; let $d'=d+m\lceil\varepsilon d\rceil$ and $g'_d$ be the exact augmented saddle at $\theta_d$. Contraction gives log-potential difference $O(\varepsilon+d^{-1})$ from the change of proportions, and the covariance and normalizer derivatives give the same orders for mean powers and log normalizers. Independently of the original paths, the added paths have average power bounded in sup norm by a constant $R'$ with probability at least $1/2$, by their bounded expected sup norm. On these bounded-power events the augmented average differs from the original by $O_R(\varepsilon+d^{-1})$, so for small $\varepsilon$ its saddle deviation is at least $\delta/2$. Lipschitz $\phi$ bounds the difference of residuals by $C_R\varepsilon d+O(1)$; the augmented tilt change has log likelihood bounded by the same order by the quadratic transform. Integrating out the added paths gives
\[
 \begin{split}
 &\E_Q[e^{s\mathcal T_d}\mathbf1_{\{\|q_d-m_g\|_\infty\ge\delta,\ \|q_d\|_\infty\le R\}}]\\
 &\qquad\le2e^{C_R\varepsilon d+O(1)}
    \E_{Q'}[e^{s\mathcal T'_{d'}}\mathbf1_{\{\|q'_{d'}-m_{g'_d}\|_\infty\ge\delta/2\}}].
 \end{split}
\]
The augmented proportions have positive limits; their negative rate for $\delta/2$ is independent of $\varepsilon$. Taking $\varepsilon$ small proves \eqref{eq:weighted-far}.

\emph{Near the saddle.} Set $Z_d=\sqrt{\theta_d}\,w_dU_d=d^{-1/2}\sum_iY_{i,d}$ with centered independent $L^2$ summands $Y_{i,d}=\sqrt{\theta_d}w_dp_i(X_i^2-C_{i,g_d}(t,t))$. Their norms have a common exponential moment by Fernique. The covariance is $K_{\theta_d}$, with uniformly vanishing trace tails because kernels and diagonal traces converge and each fixed projected trace converges.

For a finite-rank orthogonal projection $\Pi$ and unit $v$ in its range, Taylor expansion of scalar transforms to second order with the common exponential norm moment gives
\[
 \log\E_Qe^{t\langle v,Z_d\rangle}
 \le\tfrac12\langle v,K_{\theta_d}v\rangle t^2+C|t|^3/\sqrt d,\qquad
 0\le t\le c\sqrt d,
\]
using independence and $\log(1+x)\le x$. For large $d$ and small $t/\sqrt d$, the right side is at most $(q_*+\eta)t^2/2$. A finite sphere net and Chernoff then bound $Q(\|\Pi Z_d\|>u)$ by $C_\Pi\exp(-(1-\eta')^2u^2/(2(q_*+\eta)))$ for $0\le u\le c'\sqrt d$ (any $\eta,\eta'>0$, with $c'$ small enough).

Write $Y_{i,d}^\perp=(I-\Pi)Y_{i,d}$. Conditioning one independent centered term at a time gives
\[
 \E_Q\cosh(t\|(I-\Pi)Z_d\|)\le
 \exp\left(\tfrac{t^2}{2d}\sum_i
 \E_Q[\|Y_{i,d}^\perp\|^2 e^{t\|Y_{i,d}^\perp\|/\sqrt d}]\right).
\]
Indeed the Hessian of $x\mapsto\cosh(t\|x\|)$ is bounded by $t^2\cosh(t\|x\|)$ in norm, the linear Taylor term has mean zero, and $\cosh(a+b)\le\cosh(a)e^b$ for $a,b\ge0$ bounds the remainder. For $t/\sqrt d$ sufficiently small, Cauchy--Schwarz in each expectation ($E r^2e^{ar}\le(Er^2)^{1/2}(Er^2e^{2ar})^{1/2}$) and then over summands bounds the average in the exponent by $C(d^{-1}\sum_i\E_Q\|Y_{i,d}^\perp\|^2)^{1/2}$ using the common exponential norm moment. Trace tails make this smaller than any $v_*>0$ for large rank. Since $\cosh(tu)\ge e^{tu}/2$, Chernoff gives $Q(\|(I-\Pi)Z_d\|>u)\le2e^{-u^2/(2v_*)}$ for $0\le u\le c''\sqrt d$. Splitting at projected norms $(1-\eta')u,\eta'u$, first choose $\eta,\eta'$ small, then $\Pi$ large. Thus for any $\beta<1/(2q_*)$ (arbitrary finite $\beta>0$ if $q_*=0$),
\begin{equation}
 Q(\|Z_d\|>u)\le Ce^{-\beta u^2},\qquad 0\le u\le c_*\sqrt d.
 \label{eq:truncated-gaussian-tail}
\end{equation}
Choose $\beta>s/2$ using $sq_*<1$. On $\|q_d-m_g\|_\infty\le\delta$, Taylor gives $\mathcal T_d\le(1+C\delta)\|Z_d\|_2^2/2$ and $\|Z_d\|_2\le C\delta\sqrt d$. Make $\delta$ small with $a_1=s(1+C\delta)/2<\beta$ and $C\delta<c_*$. Tail integration up to $C\delta\sqrt d$ bounds the expectation by $1+\int_0^{C\delta\sqrt d}2a_1ue^{a_1u^2}Q(\|Z_d\|>u)du$, uniformly finite. The complement is exponentially negligible by \eqref{eq:weighted-far}.

\emph{Block estimates.} For a fixed balanced partition, $n_j=d/r+O(1)$ and block proportions differ from the full ones by $O(d^{-1})$. The block's own saddle at $\theta_d$ differs in log potential by $O(d^{-1})$ by contraction. By Lipschitz $\phi$, the residual difference is bounded by $C(1+\int q_j)$, as is the log density ratio of the two block tilts by the normalizer identity. Apply the full-sample bound at an exponent slightly larger than $s$ still satisfying \eqref{eq:residual-moment-range}, and H\"older. For any fixed $L$, $\E e^{L\int q_j}$ is bounded under either tilt for large blocks by the quadratic transform with coefficients divided by $n_j$ and the uniform covariance trace bound. This gives the block estimate; convexity of $\phi$ under the weights $n_j/d$ gives \eqref{eq:residual-split}.

Uniformity follows by subsequences on compact data and the crude bound for bounded dimension. For varying stationary OU rates in a compact positive interval, realize priors as rate-one stationary paths at times $\kappa_kt$ on a common bounded interval; their sup and modulus moments are uniform, and their kernels converge with the rates. Contraction and the resolvent give the same compactness and convergence of tilted kernels and traces used above.
\end{proof}

The proof also holds for a fixed positive quadrature: the grid entropy gap applies directly and the near estimate is finite-dimensional.

\subsection{A local limit in the normal coordinate}

A relative tail estimate requires the joint law within a strip of width $O(d^{-1/2})$ in $L_d$, not just a CLT. Along a sequence in Theorem~\ref{thm:sharp}, contraction and the resolvent give uniform convergence of $g_d,w_d,C_{k,g_d}$ to their limiting counterparts. Let $U$ be centered Gaussian in $L^2$ with trace-class covariance kernel $\Gamma(t,s)=2\sum_k\pi_kp_k^2 C_{k,g_\theta}(t,s)^2$. Its normal coordinate $L=\langle w_\theta,U\rangle$ has variance $s^2>0$ and density $f_L$; let $\mu_0$ be the Gaussian regression law of $U$ given $L=0$. Write $x_d=\sqrt d L_d$ and
\[
 \nu_d(dx,du)=\sqrt d\,Q((x_d,U_d)\in(dx,du)).
\]

\begin{lemma}[Joint local convergence]
\label{lem:joint-local}
On $J\times L^2([0,T])$ for every compact interval $J$, the restrictions of $\nu_d$ converge weakly as finite measures to the restriction of
\begin{equation}
 \nu(dx,du)=f_L(0)\,dx\,\mu_0(du).
 \label{eq:local-limit}
\end{equation}
The joint convergence remains valid upon adjoining $R_d$ on the left and $R_2(u)=\tfrac12\|w_\theta u\|_2^2$ on the right.
\end{lemma}

\begin{proof}
\emph{Fourier inversion.} For real $v\in L^\infty$, let $\chi_d(t;v)=\E_Q e^{\mathrm i(tL_d+\langle v,U_d\rangle)}$. The scalar triangular-array CLT gives convergence at every fixed $(t,v)$ to the characteristic function of $(L,U)$ by converging variances and bounded third moments of the unscaled summands (if the limiting variance is zero, use $L^2$ convergence). Each type-$k$ coordinate contributes the centered Gaussian quadratic transform with operator $p_k C_{k,g_d}^{1/2}M_{tw_d+v}C_{k,g_d}^{1/2}/\sqrt d$. For eigenvalues $a_j$ of a self-adjoint trace-class quadratic operator, the transform has modulus $\prod_j(1+4a_j^2)^{-1/4}$ by scalar Gaussian integration. Outside a fixed bounded interval depending on $\|v\|_\infty$, $tw_d+v$ has the sign of $t$ and absolute value at least $c_4|t|$. The lower covariance bound of Lemma~\ref{lem:central-half-estimates} gives an eigenvalue of absolute size at least $c_5|t|/\sqrt d$ for every type. Independence yields
\begin{equation}
 |\chi_d(t;v)|\le(1+c_6t^2/d)^{-d/4}.
 \label{eq:fourier-domination}
\end{equation}
For $d\ge8$ this is bounded by $(1+c_6t^2/8)^{-2}$; inside the interval use one. Thus the complex measure $\E_Q[e^{\mathrm i\langle v,U_d\rangle}\mathbf1_{\{L_d\in dl\}}]$ has a continuous density by Fourier inversion. For bounded continuous $a$ on $J$, changing variable $l=x/\sqrt d$ and using domination gives
\begin{align*}
 \int_{J\times L^2}a(x)e^{\mathrm i\langle v,u\rangle}\nu_d(dx,du)
 &=\tfrac1{2\pi}\int_Ja(x)\int_{\R}e^{-\mathrm i tx/\sqrt d}\chi_d(t;v)dt\,dx\\
 &\longrightarrow f_L(0)\int_Ja(x)dx\,\E[e^{\mathrm i\langle v,U\rangle}\mid L=0].
\end{align*}
The Gaussian regression formula gives the latter continuous weighted density at zero.

\emph{Hilbert-space tightness.} Let $\Pi_M$ project on the first $M$ of an orthonormal basis of bounded functions in $L^2$. Each half covariance is a sum of $n_{k,j}/d$ times the type kernels $2p_k^2C_{k,g_d}^2$. For the finite type list uniform kernel convergence, diagonal trace convergence and each fixed projected trace give
\[
 \sup_d\E_Q\|U_j\|_2^2<\infty,\qquad
 \lim_{M\to\infty}\sup_d\E_Q\|(I-\Pi_M)U_j\|_2^2=0 .
\]
For $A_0$ measurable in half $A$, condition on it and apply \eqref{eq:half-interval} to $L_B$ in a translated interval of length $|J|/\sqrt d$:
\begin{equation}
 \sqrt d\,Q(x_d\in J,A_0)\le C(|J|+1)Q(A_0).
 \label{eq:local-half-bound}
\end{equation}
The same holds with halves exchanged. If the total projection error exceeds $\eta$, one half exceeds $\eta/2$; thus local mass of this event is controlled by the trace tails via Markov's inequality. Similarly the total norm and local mass are bounded. To get a compact set, choose $M$'s increasing so that the respective projection tails are bounded by a sequence decreasing to zero outside sets of summably small local mass, and bound the norm. The intersection is closed and uniformly approximable by finite-dimensional bounded sets, hence compact. Tightness and the finite-coordinate Fourier limits identify the finite measures with \eqref{eq:local-limit}, since the basis separates points.

\emph{The nonlinear boundary.} Integrate \eqref{eq:local-half-bound} against $\int|U_j|^3$ and use \eqref{eq:half-moments} and $|U_A+U_B|^3\le4(|U_A|^3+|U_B|^3)$ to get
\[
 \int_{J\times L^2}\int_0^T|u(t)|^3dt\,\nu_d(dx,du)\le C_J.
\]
Thus $R_d-R_{2,d}$ tends to zero in the local finite measures by \eqref{eq:sharp-remainder}. Uniform $w_d$ convergence and tightness give joint convergence of $R_{2,d}=\tfrac12\|w_dU_d\|_2^2$ to $R_2(u)$. This result is local in $x$; the following proof controls the exponentially weighted tails.
\end{proof}

\subsection{From the local limit to the sharp probability}

\begin{proof}[Proof of Theorem~\ref{thm:sharp}]
By \eqref{eq:sharp-change-of-law}, $\sqrt d\,e^{dI}P_d(b_d)$ is the integral of $e^{\theta_d x}\mathbf1_{\{x\le R_d\}}$ under the local measures. For the unbounded $x$-tails take $s>1$ from \eqref{eq:residual-moment-range} and $r\ge3$ with $c_r:=r/(r-1)<(r-1)/(r-2)<s$, and split into $r$ balanced blocks. Their independent residuals have bounded exponential moments at $s$ by Lemma~\ref{lem:sqrt-residual}.

For $x<-M$, \eqref{eq:half-interval} bounds the mass of each unit interval uniformly under $\nu_d$; the exponential weight makes the sum vanish as $M\to\infty$ since $\theta_d\ge\theta/2$. On the positive side let $A_j$ indicate that block $j$ has minimal residual (ties broken by index), and $T_{-j}=\sum_{i\ne j}\mathcal T_i$. On $A_j$, $\theta_dR_d\le\mathcal T_j+T_{-j}\le c_r T_{-j}$. Thus $A_j\cap\{M<x_d\le R_d\}$ implies $T_{-j}\ge\theta_d M/c_r$ and $x_d\le c_r T_{-j}/\theta_d$. Given the other blocks, the normal coordinate of block $j$ (size at least $d/(2r)$) is then in a translated interval of length at most $c_r T_{-j}/(\theta_d\sqrt d)$. Dropping $A_j$ after this bound and using the block interval estimate gives
\[
 \begin{split}
 &\sqrt d\,\E_Q[e^{\theta_d x_d}\mathbf1_{A_j\cap\{M<x_d\le R_d\}}]\\
 &\qquad\le C_r\E_Q[e^{c_rT_{-j}}(1+T_{-j})
                       \mathbf1_{\{T_{-j}\ge\theta M/(2c_r)\}}].
 \end{split}
\]
Independence bounds $\E_Q e^{sT_{-j}}$ uniformly, so the sum over $j$ vanishes with $M$.

Using Lemma~\ref{lem:joint-local} on compact intervals (the boundary $x=R_2(u)$ is null under $dx\,\mu_0(du)$) and then the tail bounds gives
\[
 \lim_d\sqrt d\,e^{dI}P_d(b_d)
 =f_L(0)\int\int_{-\infty}^{R_2(u)}e^{\theta x}dx\,\mu_0(du)
 =\frac{f_L(0)}{\theta}\E[e^{\theta\|w_\theta U\|_2^2/2}\mid L=0].
\]
For $Z=w_\theta U$ of covariance $K_\theta/\theta$, spectral Gaussian integration gives mass $D_\theta=\det_F(I-K_\theta)^{-1/2}$ to the weight $e^{\theta\|Z\|^2/2}$, and weighted normalized covariance $\Sigma$. The trace bound and curvature gap ensure convergence of the integrals. The integral coordinate $L$ under that weighted law has variance $V$, hence
\begin{equation}
 f_L(0)\E[e^{\theta\|Z\|_2^2/2}\mid L=0]=\frac{D_\theta}{\sqrt{2\pi V}}.
 \label{eq:prefactor-evaluation}
\end{equation}
The conditional expectation is continuous near $L=0$ by regression into a centered Gaussian residual plus a deterministic multiple of $L$; the strict curvature gap dominates the shifted square-exponential weight there. Thus the density identity holds at zero. Uniform kernel and diagonal trace convergence, and the common gap, give convergence of actual $D,V$ via $\log D=\tfrac12\sum_{j\ge1}\operatorname{tr}(K^j)/j$ with $\operatorname{tr}(K^j)\le q_0^{j-1}\operatorname{tr}K$, and the resolvent series. This proves \eqref{eq:sharp-general}.

Insert bounded continuous functions of $(w_dU_d,R_d-x_d)$ in the same argument. In the limiting measure set $z=R_2(u)-x$; its factorization is $f_L(0)e^{\theta R_2(u)}\mu_0(du)\ e^{-\theta z}dz$ on $z\ge0$. After normalization, $z$ is exponential of rate $\theta$ independent of $w_\theta u$, which has the weighted Gaussian law conditioned on $L=0$, with regression covariance $\Sigma_\perp$. Since $R_d-x_d=d(b_d-S_d)$ this gives \eqref{eq:conditional-fluctuation}.

\emph{Moment convergence.} For the deficit on $0\le x_d\le R_d$ use $R_d-x_d\le c_rT_{-j}/\theta_d$ on $A_j$; the previous conditional bound with fourth powers has right side $C_r\E_Q[e^{c_rT_{-j}}(1+T_{-j})T_{-j}^4]$, uniformly finite by $s>c_r$. On $x_d<0$, $(R_d-x_d)^4\le C(J_A^4+J_B^4+x_d^4)$ by \eqref{eq:half-moments}. Partition into unit intervals $[-j-1,-j]$; \eqref{eq:local-half-bound} controls the half moments by their ordinary moments, and $\sum_{j\ge0}e^{-\theta j/2}(1+(j+1)^4)<\infty$.

For $\|w_dU_d\|_2^4\le C_r\sum_i\|U_i\|_2^4$ (block fluctuations normalized by $\sqrt d$), every fixed ordinary moment of $\|U_i\|_2$ is bounded by the scalar Chernoff estimates of Lemma~\ref{lem:central-half-estimates}. On the positive side, for each summand $i$, choose block $j$ of smallest residual among those other than $i$. Then $\theta_dR_d\le[(r-1)/(r-2)]T_{-j}$, and $\|U_i\|_2^4$ is measurable in the other blocks. The conditional interval bound now gives a sum over $j\ne i$ of
\[
 C_r\E_Q[\|U_i\|_2^4(1+T_{-j})e^{(r-1)T_{-j}/(r-2)}].
\]
H\"older bounds it by the ordinary moments and the $s$-exponential moment, since $(r-1)/(r-2)<s$. On the negative side, use the two-half norm bound and \eqref{eq:local-half-bound}, summing with the exponential weight. Dividing the fourth-moment bounds by the positive limiting normalizer gives uniform integrability of the linear, quadratic and mixed moments. The proof also yields the fixed positive quadrature version, with time integrals replaced by sums.
\end{proof}

\begin{corollary}[Second moment of saddle importance sampling]
\label{cor:is-moment}
In the setting of Theorem~\ref{thm:sharp}, suppose in addition that
$2\theta\rho T\max_k p_kM_k<1$.
Let $W_d=(d\Prb/dQ)\mathbf1_{\{S_d\le b_d\}}$ be the unbiased likelihood-weighted sample under the product tilt. Under this sufficient condition, define for all sufficiently large $d$, at the actual saddle,
\[
 D_{2,d}=\det_F(I-2K)^{-1/2},\qquad
 V_{2,d}=\theta_d^{-1}\langle1,K(I-2K)^{-1}1\rangle .
\]
Then
\[
 \E_Q W_d^2=
 \frac{D_{2,d}}{2\theta_d\sqrt{2\pi d V_{2,d}}}
 e^{-2dI(b_d)}(1+o(1)).
\]
In particular
\begin{equation}
 \frac{\E_Q W_d^2}{P_d(b_d)^2}
 \sim\sqrt d\,\frac{\theta\sqrt{2\pi}\,V D_2}{2D^2\sqrt{V_2}},
 \label{eq:is-relative-second-moment}
\end{equation}
where $D_2,V_2,D,V$ are evaluated at the limiting saddle. Consequently the number of independent tilted samples required for a fixed relative root-mean-square accuracy is of order $\sqrt d$.
\end{corollary}

\begin{proof}
Choose $\theta_+>\theta$ close enough that $\alpha_*>1/2$, allowing some $s>2$ in Lemma~\ref{lem:sqrt-residual}; its refined curvature bound gives $\|2K\|<1$ uniformly. By \eqref{eq:sharp-change-of-law}, $\sqrt d\,e^{2dI}\E_QW_d^2$ integrates $e^{2\theta_d x}\mathbf1_{\{x\le R_d\}}$ under $\nu_d$. The proof of Theorem~\ref{thm:sharp} applies with $r$ large enough that $2c_r<s$: the positive tail bound uses $e^{2c_rT_{-j}}$, the negative one gains decay. Local convergence and the Gaussian density argument of \eqref{eq:prefactor-evaluation}, now with covariance $\theta^{-1}K(I-2K)^{-1}$, give the limit
\[
 \frac{f_L(0)}{2\theta}\E[e^{\theta\|w_\theta U\|_2^2}\mid L=0]
 =\frac{D_2}{2\theta\sqrt{2\pi V_2}}.
\]
The common gap and trace convergence justify actual-saddle quantities as before. Division by $P_d^2$ gives \eqref{eq:is-relative-second-moment}; an average of $n$ independent samples has relative mean-square error $(\E_QW_d^2/P_d^2-1)/n$.
\end{proof}

\subsection{A relative approximation reaching the central scale}

Relative normal-tail approximations for sums in the moderate-deviation range have been studied, for instance, by Liu and Zhang \cite{LiuZhang2023} under local dependence. For the nonlinear service integral, retaining the tilted rate is useful even when the saddle tends to zero: the distance from the mean can still diverge in standard deviation units. Conversely, the leading factor in \eqref{eq:effective-centre} ceases to be accurate when that standardized distance is bounded. Let $r=\sqrt{2dI(b)}$, $u=\theta\sqrt{dV}$, and $\phi$ denote the standard normal density. Define
\begin{align}
 P_{\rm M}&=D e^{-dI(b)}e^{u^2/2}\Phi(-u),\notag\\
 P_{\rm LR}&=\Phi(-r)+\phi(r)\left(\frac{D}{u}-\frac1r\right).
 \label{eq:two-central}
\end{align}
The first expression retains the normal tail integral. The second uses the form of the Lugannani--Rice approximation \cite{LugannaniRice1980}, with the curvature factor of the nonlinear continuous event. In this section $I,D,V$ always use $b=B(\theta)$ and the actual type counts.

\begin{theorem}[Relative approximation across positive small saddles]
\label{thm:central-transition}
There exist $C<\infty$, $d_*<\infty$, and $\theta_*>0$ depending only on the fixed data, such that both approximations in \eqref{eq:two-central} are positive for $d\ge d_*$, $0<\theta\le\theta_*$, and
\begin{equation}
 \max_{P_*\in\{P_{\rm M},P_{\rm LR}\}}\left|\frac{P_d(B(\theta))}{P_*}-1\right|\le C(\theta+d^{-1/2}).
 \label{eq:central-error}
\end{equation}
Furthermore, if $\mathcal N_d$ denotes all count vectors of total size $d$ on the fixed finite type list, then
\begin{equation}
 \lim_{d\to\infty}\sup_{(n_k)\in\mathcal N_d}\sup_{0<\theta\le\theta_*}
 \max_{P_*\in\{P_{\rm M},P_{\rm LR}\}}\left|\frac{P_d(B(\theta))}{P_*}-1\right|=0.
 \label{eq:uniform-relative}
\end{equation}
In particular, the approximation remains relative even if the saddle decreases to zero while $\theta\sqrt d$ diverges.
\end{theorem}

\begin{proof}
Choose $\theta_*\le1$ such that $2C_R\theta_*<c/2$ for Lemma~\ref{lem:central-half-estimates}. First use the exact identity \eqref{eq:sharp-change-of-law}. On $L_d\le0$ its indicator is one. The additional contribution occurs on $0<\sqrt d L_d\le R_d$. Split by the larger of $J_A,J_B$. With $J_A\ge J_B$, condition on the first half and use the interval $0<\sqrt d(L_A+L_B)\le2C_RJ_A$ to bound the contribution by
\[
 \frac C{\sqrt d}\E_Q[e^{2C_R\theta J_A}(1+J_A)]\le C'd^{-1/2},
\]
using Lemma~\ref{lem:central-half-estimates} and the choice of $\theta_*$. The same bound holds with the halves exchanged. These bounds are uniform without a restriction on $\theta\sqrt d$.

Let $J(z)=e^{z^2/2}\Phi(-z)=(2\pi)^{-1/2}\int_0^\infty e^{-zv-v^2/2}dv$ for $z\ge0$. The function $x\mapsto e^{zx}\mathbf1_{\{x\le0\}}$ has total variation at most two. By integration by parts and the second bound of \eqref{eq:half-interval}, applied to $L_d/s$, we obtain
\begin{equation}
 e^{dI}P_d(B(\theta))=J(\theta\sqrt d\,s)+O(d^{-1/2}).
 \label{eq:central-j}
\end{equation}
The integral defining $J$ gives $c_7/(1+z)\le J(z)\le C_7/(1+z)$; for the lower bound integrate only over $[0,(1+z)^{-1}]$. Also $J$ is decreasing and $zJ(z)$ is increasing by a change of variables. Consequently a multiplicative $1+O(\theta)$ change of argument changes $J$ by a relative $O(\theta)$, uniformly for positive arguments.

The trace of $K$ is $O(\theta)$ uniformly by its diagonal and the bounds on $w$ and $C_{k,g}$, and $s^2\le V\le s^2/(1-\|K\|)$ by the spectral theorem. Thus $V=s^2(1+O(\theta))$ and $D=1+O(\theta)$. These comparisons and \eqref{eq:central-j} prove \eqref{eq:central-error} for $P_{\rm M}$, since $d^{-1/2}/J(\theta\sqrt d\,s)=O(d^{-1/2}+\theta)$.

To compare the second approximation, at zero tilt let $w_0=\rho/(1+\rho m_0)$ and
\[
 \sigma^2=2\sum_k\pi_{k,d}p_k^2\iint w_0(t)w_0(s)C_{k,0}(t,s)^2dt\,ds.
\]
The covariance resolvent gives uniform kernel changes of order $O(\theta)$ as $g\le\theta\rho$, hence $V=\sigma^2+O(\theta)$, and $\sigma^2$ is bounded away from zero by Lemma~\ref{lem:central-half-estimates}. Since $I(B(\theta))=\Lambda(\theta)-\theta B(\theta)$ and its derivative in $\theta$ is $\theta V(\theta)$, it follows that $I(B(\theta))=\sigma^2\theta^2/2+O(\theta^3)$, uniformly in the proportions. In particular $r/u=1+O(\theta)$ and $D/u-1/r=O(d^{-1/2})$. Rewriting
\[
 e^{dI}P_{\rm LR}=J(r)+(D/u-1/r)/\sqrt{2\pi}
\]
and using the same bounds for $J$ proves \eqref{eq:central-error} for $P_{\rm LR}$. These comparisons give positivity upon reducing $\theta_*$ and increasing $d_*$; all constants are uniform over the simplex of proportions.

For the stronger uniform relative convergence, suppose it fails and choose a violating sequence of $d$'s, count vectors and saddles. Extract a subsequence with convergent proportions and $\theta_d\to\theta\in[0,\theta_*]$. If $\theta=0$, \eqref{eq:central-error} contradicts the violation. If $\theta>0$, Theorem~\ref{thm:sharp} applies with actual counts. Here $u$ and $r$ both grow as $\sqrt d$, and the normal integral satisfies $J(z)=(2\pi)^{-1/2}(z^{-1}+O(z^{-3}))$ as $z\to\infty$, by the integral and $|1-e^{-v^2/2}|\le v^2/2$. Substitution into \eqref{eq:two-central} gives $P_{\rm M}/P_{0,d}\to1$ and $P_{\rm LR}/P_{0,d}\to1$ (for the latter, the $1/r$ terms cancel). This again contradicts the violation and proves \eqref{eq:uniform-relative}.
\end{proof}

The restriction to a Gaussian tail factor at the central scale is substantive. For example, in one type of unit-variance stationary OU paths with $T=\rho=p=\kappa=1$, $V_0=(1+e^{-2})/4$ at zero tilt. If $\theta_d=(V_0d)^{-1/2}$, the expansions just derived give $r,u\to1$ and $D\to1$. By Theorem~\ref{thm:central-transition}, $P_d(B(\theta_d))\to\Phi(-1)$, whereas $P_{0,d}\to\phi(1)$. Thus a relative error tending to zero for the leading boundary factor alone cannot hold uniformly down to zero saddle. The Gaussian tail in \eqref{eq:two-central} accounts for this transition.

\begin{corollary}[Conditional fluctuations at small saddles]
\label{cor:moderate-conditioning}
Suppose the proportions converge and $\theta_d\to0$. Use the actual saddle and $b_d=B(\theta_d)$ as above. At the limiting proportions let
\[
 \Sigma_0(t,s)=2w_0(t)w_0(s)\sum_k\pi_kp_k^2 C_{k,0}(t,s)^2,
 \qquad V_0=\langle1,\Sigma_0 1\rangle>0,
\]
where $w_0=\rho/(1+\rho m_0)$. Let $Z_{0,\perp}$ be centered Gaussian with covariance
$\Sigma_0-(\Sigma_0 1)\otimes(\Sigma_0 1)/V_0$.
Under $S_d\le b_d$, the following joint limits hold.
\begin{enumerate}
\item If $\theta_d\sqrt d\to\infty$, then
\[
 \big(w_d U_d,\,\theta_d d(b_d-S_d)\big)
 \Longrightarrow (Z_{0,\perp},E_1).
\]
Here $E_1$ is an independent exponential of rate one.
\item If $\theta_d\sqrt d\to c\in[0,\infty)$, then
\[
 \big(w_dU_d,\,\sqrt d(b_d-S_d)\big)
 \Longrightarrow
 \left(Z_{0,\perp}-\frac{\Sigma_0 1}{V_0}Y_c,\,Y_c\right),
\]
where $Y_c$ is independent of $Z_{0,\perp}$ and has density on $y\ge0$
given by
\[
 \frac{\exp(-cy-y^2/(2V_0))}{\sqrt{2\pi V_0}\,J(c\sqrt{V_0})}.
\]
\end{enumerate}
\end{corollary}

\begin{proof}
\emph{Part (i).} Put $a_d=\theta_d\sqrt d\to\infty$ and $z_d=-a_d L_d$. Under $Q$ consider the local measures $a_d Q((z_d,U_d)\in(dz,du))$. The proof of Lemma~\ref{lem:joint-local} applies with limiting covariance at zero tilt: Fourier inversion now has the factor $e^{\mathrm i t z/a_d}$, which tends to one on each compact interval, while the dominating bound \eqref{eq:fourier-domination} is unchanged. The half-sample bound for a compact interval $J$ is now $a_d Q(z_d\in J,A_0)\le C(|J|+a_d/\sqrt d)Q(A_0)\le C(|J|+1)Q(A_0)$, by \eqref{eq:half-interval}. This gives the same Hilbert-space tightness. The local measures therefore converge on compact intervals to $f_{L_0}(0)\,dz\,\mu_{0,0}(du)$, where $L_0$ is the normal coordinate at zero tilt and $\mu_{0,0}$ the corresponding Gaussian profile conditioned on $L_0=0$. By \eqref{eq:half-moments} and the modified half-sample bound, the local integral of $R_d$ on compact intervals is bounded. Thus $\theta_dR_d\to0$ in these local finite measures.

The change of law is $a_d e^{dI}P_d(b_d)=a_d\E_Q[e^{-z_d}\mathbf1_{\{z_d\ge-\theta_dR_d\}}]$. The part with negative $z_d$ is the positive-boundary contribution in the proof of Theorem~\ref{thm:central-transition}, bounded here by $a_d O(d^{-1/2})=O(\theta_d)\to0$. For $z_d\ge0$ the indicator is one, and the integral over $z_d>M$ is bounded by a constant times $e^{-M}$ by the interval estimate and a sum over unit intervals. On compact positive intervals the local limit applies, also with $\theta_dR_d$ adjoined. It follows that the normalized weighted measure converges to the product of $e^{-z}dz$ on $z\ge0$ and $\mu_{0,0}$. Since the scaled service deficit is $z_d+\theta_dR_d$ and $w_d\to w_0$, Gaussian regression gives the stated limit. The curvature weighting tends to one in this regime, leaving the unweighted projected covariance.

\emph{Part (ii).} If $a_d\to c<\infty$, the scalar triangular-array central limit theorem in the proof of Lemma~\ref{lem:joint-local}, together with its (unweighted) covariance trace-tail bounds, gives $w_dU_d\Rightarrow Z_0$ under $Q$, with $Z_0$ centered Gaussian of covariance $\Sigma_0$. Consequently $L_d\Rightarrow L_0=\langle1,Z_0\rangle\sim N(0,V_0)$. Lemma~\ref{lem:central-half-estimates} gives $\E_Q R_d\le C$, so $R_d/\sqrt d\to0$ in probability. In the likelihood identity \eqref{eq:sharp-change-of-law}, the contribution of $0<\sqrt d L_d\le R_d$ is $O(d^{-1/2})$ by the positive-side estimate in Theorem~\ref{thm:central-transition}. The same bound holds with a bounded test function of the joint variables. On $L_d\le0$ the weight $e^{a_dL_d}$ is bounded by one; hence joint weak convergence gives the limiting law of $(Z_0,-L_0)$ weighted by $e^{cL_0}\mathbf1_{\{L_0\le0\}}$. The boundary $L_0=0$ has probability zero, and its total weight is $J(c\sqrt{V_0})>0$. Gaussian regression writes $Z_0=Z_{0,\perp}+(\Sigma_0 1)L_0/V_0$ with independent $Z_{0,\perp}$ and $L_0$. Only the law of $L_0$ changes under the weight. Taking $Y_c=-L_0$ under this normalized law gives the stated density and joint limit.
\end{proof}

\subsection{The first correction at the central scale}
\label{sec:edgeworth}

The relative limit in Theorem~\ref{thm:central-transition} does not distinguish the two normal-tail formulas at the level of their first correction. For fixed rational type proportions the Lugannani--Rice form also matches that correction. We include the coefficients explicitly because the typical service is not the exact finite-sample mean.

Use the prior kernels $C_k=C_{k,0}$, put $w=w_0=\rho/(1+\rho m_0)$, and define
\begin{gather*}
 K_0^*(s,t)=2\sum_k\pi_kp_k^2C_k(s,t)^2,\qquad h(t)=\int_0^T K_0^*(t,s)w(s)ds,\\
 \sigma^2=\int_0^T w(t)h(t)dt,\qquad \mu=-\tfrac12\int_0^T w(t)^2K_0^*(t,t)dt,\\
 \begin{aligned}
 \gamma={}&8\sum_k\pi_kp_k^3\int_{[0,T]^3}w(s)w(t)w(v)\\
 &\qquad\times C_k(s,t)C_k(t,v)C_k(v,s)\,ds\,dt\,dv-3\int_0^T w(t)^2h(t)^2dt.
 \end{aligned}
\end{gather*}

\begin{proposition}[Central distribution correction]
\label{prop:edgeworth}
Fix rational proportions and let $d$ run through their integer multiples. Then, uniformly over $x\in\mathbb R$,
\begin{equation}
 P_d\left(b_0+\frac{\sigma x}{\sqrt d}\right)=\Phi(x)+\frac{\phi(x)}{\sqrt d}
     \left[\frac{\gamma}{6\sigma^3}(1-x^2)-\frac\mu\sigma\right]+o(d^{-1/2}).
 \label{eq:edgeworth}
\end{equation}
On each bounded interval $-H\le x\le0$, $P_{\rm LR}$ in \eqref{eq:two-central} has the same expansion, using its continuous extension at $b=b_0$. The mean satisfies $\E S_d=b_0+\mu/d+O(d^{-3/2})$.
\end{proposition}

\begin{proof}
First expand the saddle at zero. In the uniform kernel norm, the covariance resolvent has the expansion
\[
 C_{k,g}=C_k-2p_k C_kM_gC_k+4p_k^2 C_kM_gC_kM_gC_k+O(\|g\|_\infty^3).
\]
The Neumann series is locally convergent, and the fixed-point equation is analytic at $(\theta,g)=(0,0)$ with identity derivative in $g$ after moving its right side to the left. In particular $g=\theta w+\theta^2w^2h+O(\theta^3)$. The first diagonal covariance change gives $m_g-m_0=-\theta h+O(\theta^2)$. At quadratic order, insertion of the potential correction into the linear covariance term contributes $-\int w^2h^2$ to $B$, the quadratic covariance term contributes the triple-integral term of $\gamma/2$, and the quadratic term of the service logarithm contributes $-\tfrac12\int w^2h^2$. This gives, with differentiable local remainders,
\begin{align*}
 B(\theta)&=b_0-\sigma^2\theta+\gamma\theta^2/2+O(\theta^3),& V&=\sigma^2-\gamma\theta+O(\theta^2),\\
 I(B(\theta))&=\sigma^2\theta^2/2-\gamma\theta^3/3+O(\theta^4),&D&=1-\mu\theta+O(\theta^2).
\end{align*}
Here we used $V=-B'$, $I(B(\theta))=\int_0^\theta tV(t)dt$ and $\log D=\tfrac12\operatorname{tr}K+O(\theta^2)$.

For the distribution expansion write $\pi_k=l_k/r_0$ with integer $l_k$, $r_0=\sum_k l_k$, and $d=nr_0$. Regroup the paths into $n$ iid blocks of composition $(l_k)_k$. Let $Z_i(t)$ be the centered block power, normalized by $r_0$, so $q_d=m_0+\overline Z$ with $\overline Z=n^{-1}\sum_i Z_i$. Define
\[
 L_i=\int wZ_i,\qquad H_{ij}=-\int w^2Z_iZ_j,\qquad a=\tfrac12\E H_{11}.
\]
For $i\ne j$ the kernel $H_{ij}$ is degenerate, since $\E(H_{ij}\mid Z_i)=0$. Taylor expansion yields $S_d=b_0+U_n+a/n+\mathcal R_n$, where $U_n$ is the degree-two $U$-statistic with symmetric kernel $(L_i+L_j+H_{ij})/2$, normalized by $\binom n2^{-1}$, and $\|\mathcal R_n\|_2=O(n^{-3/2})$. To verify the remainder, the third derivative of the service logarithm is bounded on nonnegative powers, so the integrated cubic error is at most $C\int |\overline Z|^3$. Its $L^2$ norm is $O(n^{-3/2})$, since in the expansion of the sixth moment only index patterns in which every independent centered summand occurs at least twice survive, giving $\E|\overline Z(t)|^6=O(n^{-3})$ uniformly in time. The centered diagonal terms $\sum_i(H_{ii}-\E H_{ii})/(2n^2)$ have this same order; the off-diagonal coefficient difference $1/n^2-1/[n(n-1)]$ is $O(n^{-3})$, and the norm of the sum of the degenerate pair terms is $O(n)$ because distinct unordered pairs are uncorrelated.

The projection of this $U$-statistic is $L/2$ with positive variance. All moments needed for the one-term expansion in Jing and Wang \cite[Theorem~1.1]{JingWang2003} hold: the projection has a finite third moment and the kernel a finite $5/3$ moment. Moreover $L$ is a nonzero centered Gaussian quadratic form (a sum over the independent paths in a block), with at least one strictly positive eigenvalue in its covariance representation. Its characteristic function tends to zero in modulus at infinity by the Gaussian quadratic integral, so its distribution is nonlattice. The block path space is standard Borel and can be represented by real observations for that theorem. Put $s^2=\E L^2$. In its normalization the third coefficient is $(\E L^3+3\E[L_1L_2H_{12}])/s^3$. Wick's identity gives
\[
 s^2=\sigma^2/r_0,\quad a=\mu/r_0,\quad
 \E L^3+3\E[L_1L_2H_{12}]=\gamma/r_0^2.
\]
Indeed the third centered moment of three Gaussian squares is $8C(s,t)C(t,v)C(v,s)$, and $\E[LZ(t)]=h(t)/r_0$. The uniform expansion for $\sqrt nU_n/s$ thus transfers to $\sqrt d(S_d-b_0)/\sigma$: by Chebyshev $\Prb(|\sqrt n\mathcal R_n/s|>n^{-2/3})=O(n^{-2/3})$, so bounding the CDF between shifts of size $n^{-2/3}$ changes its smooth normal and correction terms by $o(n^{-1/2})$. The deterministic shift $a/n=\mu/d$ contributes $-\mu\phi(x)/(\sigma\sqrt d)$. This proves \eqref{eq:edgeworth}; expectations of the same expansion give the mean formula.

Finally write $b=b_0-\sigma z/\sqrt d$, $0<z\le H$. Series division using the saddle expansions gives
\[
 r=z+\frac{\gamma z^2}{6\sigma^3\sqrt d}+O(z^3/d),\qquad
 D/u-1/r=\frac1{\sqrt d}\left(-\frac\mu\sigma+\frac{\gamma}{6\sigma^3}\right)+O(z/d).
\]
The singular $1/\theta$ terms cancel before letting $\theta$ decrease to zero, so the formula extends continuously at $z=0$. Inserting these expansions in $P_{\rm LR}$ gives the right side of \eqref{eq:edgeworth} at $x=-z$ uniformly on $[0,H]$.
\end{proof}

The continuous extension at $b_0$ is an asymptotic approximation and need not belong to $[0,1]$ at small counts.

\section{A quantitative relative error for Ornstein--Uhlenbeck paths}
\label{sec:ou-quantitative}

For OU paths the exact radial likelihood permits a relative error of order $d^{-1}$: after normal-coordinate inversion, the first correction, including the curved boundary and Jacobian, is odd under a centered Gaussian law. Related symmetry cancellations occur in moderate deviations of finite-dimensional quadratic norms \cite{FangLiuShao2023}; here the saddle is positive and the expansion uses the full radial diffusion.

\begin{theorem}[Uniform relative order for OU types]
\label{thm:ou-relative}
Fix $T,\rho>0$ and a finite number $m$ of types with $C_{k,0}(s,t)=e^{-\kappa_k|t-s|}$. Consider any compact family of positive powers and decay rates, actual proportions $\pi_k\ge\gamma>0$ and multipliers $0<\theta_-\le\theta\le\theta_+<\infty$. For constants $C<\infty$ and $d_0$ depending only on these data, all integer count vectors in the family and $d\ge d_0$ satisfy
\begin{equation}
 \left|\frac{P_d(B(\theta))}{P_{0,d}(B(\theta))}-1\right|\le C/d.
 \label{eq:ou-relative}
\end{equation}
The same relative order holds for both approximations in \eqref{eq:two-central} (with possibly a different $C,d_0$).
\end{theorem}

\subsection{Exact radial representation}

For each type write $c_k(t)=C_{k,g}(t,t)$. The tilted law is a linear diffusion. Indeed, for the prior stationary OU process of rate $\kappa_k$, the backward expectation of the remaining quadratic weight has the form
\begin{gather*}
 h_k(t,x)=\exp(A_k(t)-r_k(t)x^2/2),\\
 r_k'=2\kappa_k(r_k+r_k^2)-2p_kg(t),\qquad A_k'=\kappa_k r_k,
 \quad r_k(T)=A_k(T)=0.
\end{gather*}
In reversed time $r_k$ is nonnegative and bounded by the vector-field signs at zero and a sufficiently large constant. Substitution in the backward equation and It\^o's formula with the bounded exponential weight verify $h_k$. The change of law has initial variance $(1+r_k(0))^{-1}$, drift $-\kappa_k(1+2r_k)X$ (by the stochastic product rule) and variance equation $c_k'=2\kappa_k-2\kappa_k(1+2r_k)c_k$.

\medskip\noindent\textit{Relation to the radial path action.}
For a unit-variance stationary OU process of rate $\kappa$, the average of $n$ independent squared coordinates solves
\[
 dY=2\kappa(1-Y)dt+\sqrt{8\kappa Y/n}\,dW,\qquad nY(0)\sim\chi_n^2.
\]
For a prescribed positive normalized initial value, the dynamic action per coordinate in Donati-Martin et al.~\cite[Proposition~1.2 and Section~6.1]{DonatiMartinRouaultYorZani2004}, after a time change, is the integral in
\begin{equation}
 J_\kappa(c)=\frac{c(0)-1-\log c(0)}2
   +\int_0^T\frac{[\dot c(t)-2\kappa(1-c(t))]^2}{16\kappa c(t)}\,dt.
 \label{eq:stationary-radial-action}
\end{equation}
The initial term is the Cram\'er rate of $\chi_n^2/n$. Use this action on nonnegative absolutely continuous paths, infinity otherwise, with infinite initial cost at zero and lower-semicontinuous extension of the integrand there. Finite-action paths are bounded away from zero: before a first zero $t_*>0$, expansion of the square gives
\[
 \int_0^t\frac{[\dot c-2\kappa(1-c)]^2}{16\kappa c}\,ds
 \ge -\tfrac14[\log c(t)-\log c(0)-c(t)+c(0)].
\]
The right side diverges at $t_*$. Hence $c>0$ on the compact interval and $\dot c\in L^2$. For such $c$, let $Q^c$ be the linear Gaussian diffusion of initial variance $c(0)$, noise $\sqrt{2\kappa}$ and drift $a(t)X$, $a=(\dot c-2\kappa)/(2c)\in L^2$, so its variance is $c$. After changing the initial law, its drift entropy relative to stationary OU is
\[
 \frac1{4\kappa}\E_{Q^c}\int_0^T(a(t)+\kappa)^2X(t)^2dt
 =\int_0^T\frac{[\dot c-2\kappa(1-c)]^2}{16\kappa c}\,dt.
\]
To justify Girsanov, stop at $|X|=L$; the integrated squared drift difference is bounded by $L^2\int(a+\kappa)^2$. Both linear diffusions are nonexplosive by integrability of their deterministic drifts, so the stopped changes exhaust the paths. Under $Q^c$, the log drift likelihood is the quadratic term before expectation plus $\int(a+\kappa)X/\sqrt{2\kappa}\,dW_{Q^c}$ of mean zero by square integrability. Adding initial entropy gives $D(Q^c\Vert P)=J_\kappa(c)$. Thus Theorem~\ref{thm:entropy-stability} yields
\[
 I_H(b)=\min_{\substack{c_k\ge0,\ c_k\text{ absolutely continuous}\\
     \int_0^T\log(1+\rho\sum_k\pi_kp_k c_k)\,dt\le b}}
     \sum_k\pi_k J_{\kappa_k}(c_k).
\]
Indeed each finite-action profile is realized by this law, and the tilted variances attain the minimum uniquely (their dynamic integrand is $\kappa_k r_k^2c_k$). This solves the stationary radial minimization and identifies the individual conditioned laws.

Returning to the tilted paths, division by $\sqrt{c_k(t)}$ and clock $\tau_k(t)=\int_0^t\kappa_k/c_k(s)ds$ give stationary variance-one OU of rate one. The resolvent gives a common positive lower bound for $c_k$ and hence bounded terminal clocks.

In clock $\tau_k$, let $z_k$ be the radius of the $n_k$ standardized coordinates minus $\sqrt{n_k}$. For $n_k\ge2$ it stays in $(-\sqrt{n_k},\infty)$: a deterministic scaling and time change turn the original vector OU into Brownian motion in that dimension, for which a point is polar. Its diffusion coefficient is $\sqrt2$ and its drift is
\[
 b_n(z)=\frac{n-1}{\sqrt n+z}-(\sqrt n+z)=-2z+\delta_n(z),\qquad
 \delta_n(z)=\frac{z^2-1}{\sqrt n+z}.
\]
Use the reference law $\mathcal G_0$ of independent stationary OU paths $Z_k$ of rate two and variance $1/2$ in their respective clocks (covariance $\tfrac12e^{-2|\tau_k(t)-\tau_k(s)|}$ in physical time). For $N=\sqrt n$, set
\begin{gather*}
 l_n(z)=c_n+(n-1)\log(1+z/N)-Nz+z^2/2,\\
 c_n=\tfrac12\log\pi+(1-n/2)\log2-\log\Gamma(n/2)+\tfrac{n-1}2\log n-n/2.
\end{gather*}
The initial shifted chi density relative to the reference density is $e^{l_n(z)}$, and $l_n'=\delta_n$. Girsanov's formula stopped on compact subintervals of the physical domain and It\^o's formula give the exact path log likelihood
\begin{equation}
 \begin{gathered}
 \mathscr B(Z)=\sum_k\left[\frac{l_{n_k}(Z_k(0))+l_{n_k}(Z_k(T))}{2}
     +\int_0^{\tau_k(T)} k_{n_k}(Z_k)\,d\tau\right],\\
 k_n(z)=z\delta_n(z)-\frac{\delta_n'(z)}2-\frac{\delta_n(z)^2}4.
 \end{gathered}
 \label{eq:ou-radial-likelihood}
\end{equation}
The integral is in the type clock. Before It\^o's rule the stopped log likelihood is $l_n(Z_0)+\tfrac12\int\delta_n(Z)dZ+\int(Z\delta_n-\delta_n^2/4)d\tau$. Restrict the initial value to $[-\sqrt n+1/L,L]$ and stop at its first exit: the bounded drift difference permits the stopped change of law \cite{RevuzYor1999}. Integrate against the initial density ratio and let $L\to\infty$. The no-exit events exhaust physical paths since the radius stays positive and finite; monotone convergence gives $dQ/d\mathcal G_0=e^{\mathscr B}$ on physical radius paths and zero otherwise (here $Q$ denotes the radial image of the product tilt).

Put $\varepsilon=d^{-1/2}$, $w_k=\rho p_kc_k/(1+\rho m_g)$, $f_k=2\sqrt{\pi_k}w_k$, and write the linear and quadratic radial contributions as
\[
 B_1=\sum_k f_kZ_k,\qquad B_2=\sum_kw_kZ_k^2,\qquad L(Z)=\int_0^T B_1(t)dt.
\]
The power ratio is $z=(1+\rho q_d)/(1+\rho m_g)=1+\varepsilon B_1+\varepsilon^2 B_2$. The exact scaled residual is $\mathcal T=\theta R_d=d\theta\int\phi(z)$, with quadratic reference term $\mathcal T_2=\theta\int B_1^2/2$ (the multiplier is included here, unlike $R_{2,d}$ in \eqref{eq:sharp-remainder}). Under $\mathcal G_0$, $B_1$ has covariance $K/\theta$ because $C_{k,g}(s,t)=\sqrt{c_k(s)c_k(t)}e^{-|\tau_k(t)-\tau_k(s)|}$. Thus the weighted Gaussian law $\mathcal G$ has density $e^{\mathcal T_2}/D$ relative to $\mathcal G_0$, and $B_1$ under $\mathcal G$ has covariance $\theta^{-1}K(I-K)^{-1}$. Under $\mathcal G$ regress
\[
 x=L(Z)\sim N(0,V),\qquad Z=\xi+\beta x,\qquad
 \beta_k(t)=\operatorname{Cov}_{\mathcal G}(Z_k(t),L)/V,\qquad L(\xi)=0,
\]
with $\xi$ of law $\nu$ independent of $x$. All $\beta_k$ are positive: weighting multiplies the positive reference cross-covariance with $B_1$ against $(I-K)^{-1}1>0$. The bounds on $V$, the reference sup-norm exponential square moment on bounded clocks and the curvature gap (giving $\E_{\mathcal G_0}e^{s_0\mathcal T_2}\le C$ for some $s_0>1$) give by H\"older a uniform exponential square moment under $\mathcal G$. Cauchy--Schwarz then bounds $\beta$ uniformly, and $|x|\le C\|Z\|_\infty$ implies
\begin{equation}
 \sup\E_\nu e^{c_0\|\xi\|_\infty^2}<\infty
 \label{eq:ou-residual-moments}
\end{equation}
for some $c_0>0$, with norm the maximum over type paths. All bounds are uniform on the compact family.

\subsection{Localization and cancellation of the first correction}

For fixed $\xi$ all physical $x$ form an interval since $\beta_k(t)>0$. For the reconstructed radius let $Y(x)=\sqrt d(S_d-b)$. Its derivative is
\[
 J=\int_0^T\frac{H+\varepsilon D_1}{z},\qquad H=\sum_k f_k\beta_k,\quad
 D_1=2\sum_k w_k\beta_k Z_k,\quad \int_0^T H=1 .
\]
It is positive on that interval because $H+\varepsilon D_1=2\sum_kw_k\beta_k(\sqrt{\pi_k}+\varepsilon Z_k)>0$. Put $F=\mathcal T-\mathcal T_2+\mathscr B-\log J$. In \eqref{eq:sharp-change-of-law}, write $u=d(S_d-b)=Y/\varepsilon$, use $dQ/d\mathcal G_0=e^{\mathscr B}$ on physical paths and $d\mathcal G/d\mathcal G_0=e^{\mathcal T_2}/D$, then substitute $x$ in the normal density under $\mathcal G$:
\begin{equation}
 \frac{P_d(b)}{P_{0,d}(b)}=
 \theta\int_{-\infty}^0 e^{\theta u}\E_\nu\big[\mathbf1_{\{\text{physical root exists}\}}
       e^{F(\xi+\beta x)-x^2/(2V)}\big]du. 
 \label{eq:ou-boundary-integral}
\end{equation}
Here $x$ solves $Y(x)=\varepsilon u$ uniquely on the physical interval, if it exists. The normalization $P_{0,d}$ uses precisely the Gaussian mass $D$, normal density at zero and factor $\varepsilon/\theta$.

\emph{Localization under the tilted law.} Restrict \eqref{eq:ou-boundary-integral} to $\|\xi\|_\infty\le M=A\sqrt{\log d}$ and $-U\le u\le0$, $U=B_0\log d$, with large constants. Under $Q$, a shifted radius has exit probability from $[-M,M]$ at most $Ce^{-cM^2}$ for $M\le\min_k\sqrt{n_k}/2$. Its initial chi norm is 1-Lipschitz in the standard normals, with mean within one of $\sqrt n$ (variance at most one and second moment $n$), giving a common shifted exponential square moment. On $|z|\le\sqrt n/2$,
\[
 zb_n(z)=-\frac{2\sqrt n+z}{\sqrt n+z}z^2-\frac{z}{\sqrt n+z}
 \le-\tfrac32 z^2+C.
\]
For small fixed $\alpha>0$, the generator of $e^{\alpha z^2}$ is at most $Ce^{\alpha z^2}$. The resulting exponential supermartingale stopped at the first exit bounds the exit probability over a bounded clock interval by $Ce^{-\alpha M^2}$, including initial exits; sum over types.

Lemma~\ref{lem:sqrt-residual} gives $\E_Q e^{s\mathcal T}\le C$ for some $s>1$, uniformly. The ratio in \eqref{eq:ou-boundary-integral} is a bounded constant times $\sqrt d\,\E_Q[e^{\theta u+\mathcal T}\mathbf1_{\{u\le0\}}]$. Since $\|\xi\|_\infty=\|Z-\beta L(Z)\|_\infty\le C\|Z\|_\infty$, H\"older and the true exit estimate bound the part outside the box by $C\sqrt d e^{-c'M^2}$; the part $u<-U$ is at most $C\sqrt d e^{-\theta_-U}$. Both are $O(d^{-2})$ with sufficiently large $A,B_0$.

\emph{Normal-coordinate inversion.} On $\|\xi\|_\infty\le M$, expansions on the tube $\|Z\|_\infty\le2M$ give $J=1+O(\varepsilon M)$ and $|Y(0)|\le C\varepsilon\|\xi\|_\infty^2$ since $\int B_1(\xi)=0$. For $|u|\le U$, evaluate $Y$ at $\pm2\varepsilon(|u|+C\|\xi\|_\infty^2)$: $J\ge1/2$ brackets the desired value. These points remain in the tube since $\varepsilon\log d=o(M)$ and $\beta$ is bounded; the tube is physical by $n_k\ge\gamma d$. The unique root satisfies
\begin{equation}
 |x|\le C\varepsilon(|u|+\|\xi\|_\infty^2).
 \label{eq:ou-root-bound}
\end{equation}

\emph{Expansion and cancellation.} Stirling's formula with remainder yields $c_n=O(n^{-1})$. The Taylor series for $\log(1+z/\sqrt n)$ and the rational expression for $\delta_n$ give, uniformly on the tube $|z|\le2M$,
\begin{align*}
 l_n(z)&=n^{-1/2}(z^3/3-z)+O(n^{-1}(1+|z|^4)),\\
 k_n(z)&=n^{-1/2}(z^3-2z)+O(n^{-1}(1+|z|^4)).
\end{align*}
For the latter, expand $\delta_n=(z^2-1)/\sqrt n+O(n^{-1}(|z|^3+|z|))$, $\delta_n'=2z/\sqrt n+O(n^{-1}(1+z^2))$ and use $\delta_n^2=O(n^{-1}(1+z^4))$. The remainders are uniform since $|z|/\sqrt n\le1/2$, and $n_k^{-1/2}=\varepsilon/\sqrt{\pi_k}$ with $\pi_k^{-1/2}\le\gamma^{-1/2}$. Expanding $\phi(1+v)=v^2/2-v^3/3+O(v^4)$ at $v=\varepsilon B_1+\varepsilon^2 B_2$ and the Jacobian $(H+\varepsilon D_1)/z$ using $\int H=1$ gives
\begin{align*}
 \mathcal T-\mathcal T_2&=\varepsilon\theta\int_0^T(B_1B_2-B_1^3/3)+O(\varepsilon^2(1+\|Z\|_\infty^4)),\\
 -\log J&=-\varepsilon\int_0^T(D_1-HB_1)+O(\varepsilon^2(1+\|Z\|_\infty^2)).
\end{align*}
Consequently $F(Z)=\varepsilon F_1(Z)+O(\varepsilon^2(1+\|Z\|_\infty^4))$, where the full first-order term is
\begin{equation}
 \begin{aligned}
 F_1(Z)={}&\sum_k\frac1{\sqrt{\pi_k}}\left[
     \tfrac12\sum_{t\in\{0,T\}}(Z_k(t)^3/3-Z_k(t))
     +\int_0^{\tau_k(T)}(Z_k^3-2Z_k)\,d\tau\right]\\
 &+\theta\int_0^T(B_1B_2-B_1^3/3)\,dt-\int_0^T(D_1-HB_1)\,dt .
 \end{aligned}
 \label{eq:ou-first-correction}
\end{equation}
This odd polynomial functional has degree at most three and bounded coefficients and integration lengths. Its change along $\beta x$ is bounded by $C|x|(1+\|\xi\|_\infty^2+|x|^2)$. Using \eqref{eq:ou-root-bound} and the lower bound for $V$ gives at the root
\[
 F(\xi+\beta x)-x^2/(2V)=\varepsilon F_1(\xi)+\varepsilon^2\mathcal E_d(\xi,u),
 \qquad |\mathcal E_d(\xi,u)|\le C(1+\|\xi\|_\infty^6+|u|^4).
\]
Indeed the cubic change is bounded by $C\varepsilon^2(|u|+\|\xi\|^2)(1+\|\xi\|^2+|x|^2)$ and $x^2$ by $C\varepsilon^2(|u|^2+\|\xi\|^4)$. The exponent is uniformly $o(1)$ since $M=O(\sqrt{\log d})$, $U=O(\log d)$. Hence
\[
 e^{F(\xi+\beta x)-x^2/(2V)}=1+\varepsilon F_1(\xi)
       +O(\varepsilon^2(1+\|\xi\|_\infty^{12}+|u|^8)).
\]
For every parameter choice, $\nu$ is centered Gaussian and the box $\|\xi\|_\infty\le M$ is symmetric. Thus the integral of $F_1(\xi)$ is exactly zero. The polynomial remainder is uniformly integrable by \eqref{eq:ou-residual-moments} and $\theta e^{\theta u}$ on $u\le0$. The omitted mass of the constant term is $O(d^{-2})$ by the same Gaussian and exponential tails with large $A,B_0$. Together with the discarded true mass, \eqref{eq:ou-boundary-integral} gives \eqref{eq:ou-relative}.

On the compact separated family, $I,V$ are bounded below and $D,V$ above, so $r,u$ in \eqref{eq:two-central} grow uniformly as $\sqrt d$. The normal-tail expansion used in Theorem~\ref{thm:central-transition} gives $P_{\rm M}/P_{0,d}=1+O(d^{-1})$ and, by cancellation of $1/r$, the same for $P_{\rm LR}$. This completes the proof. Uniformly separated internal thresholds equivalently give compact positive multiplier bounds by Proposition~\ref{prop:global-response}.

\section{Reliable service quantiles}
\label{sec:quantiles}

The probability result can be inverted at an exponentially small target, which is useful when the accumulated service is to be guaranteed at a prescribed reliability level. We state the inverse correction for the actual type proportions.

For every finite $d$, the service law is continuous and has a strictly increasing CDF on $(0,\infty)$. To verify this, select one of its Gaussian coordinates and regress its whole path on a nondegenerate evaluation, writing $X(t)=a f(t)+Y(t)$, where $a$ is standard normal, $f(t)$ is continuous and nonzero on some interval, and $Y$ is independent of $a$. Condition on $Y$ and all other paths. The service is a real analytic function of $a$: on any compact real range the arguments of the logarithms are at least one and their polynomial coefficients are bounded, so the logarithm has a uniform local analytic continuation. The function tends to infinity as $|a|\to\infty$ on account of the interval where $|f|$ stays positive. Each level set of this nonconstant real analytic function is discrete, hence conditional probability of a given level is zero. For strict increase, the centered residual and other independent Gaussian paths have positive probabilities of arbitrary uniform neighborhoods of zero. Indeed, by separability some ball of any prescribed radius has positive probability; the difference of two independent draws in that ball, divided by $\sqrt2$, has the original centered Gaussian law and lies near zero. When these paths vanish, varying $a$ in $[0,\infty)$ makes the service vary continuously from zero to infinity. Continuity in all paths and a positive probability of an appropriate neighborhood of $a$ then give positive probability to every open subinterval of $(0,\infty)$.

\begin{proposition}[The first quantile correction]
\label{prop:quantile}
Let the proportions converge as in Theorem~\ref{thm:continuous-exponent}. Fix $r>0$ and let $b_r$ be the unique internal solution of $I_H(b_r)=r$, which exists by Proposition~\ref{prop:global-response}. Let $b_{r,d}$ solve $I_{H_d}(b_{r,d})=r$ for the actual counts; write $\theta_{r,d},D_{r,d},V_{r,d}$ for its saddle quantities and $A_{r,d}=D_{r,d}/(\theta_{r,d}\sqrt{2\pi V_{r,d}})$. The unique quantile $b_d^*$ satisfying $P_d(b_d^*)=e^{-dr}$ obeys
\begin{equation}
 b_d^*=b_{r,d}+\frac{\tfrac12\log d-\log A_{r,d}}{d\theta_{r,d}}+o(d^{-1}).
 \label{eq:quantile-expansion}
\end{equation}
Under the OU relative-order assumptions the remainder is $O((\log d)^2/d^2)$ uniformly on compact positive ranges of reliability levels.
\end{proposition}

\begin{proof}
The saddle and rate depend continuously on the proportions, so $b_{r,d}\to b_r$ by the strict derivative $I'_H=-\theta$. On a fixed small neighborhood of $b_r$ all limiting multipliers are bounded away from zero and infinity. Theorem~\ref{thm:sharp}, by a compact-subsequence argument, implies the relative approximation uniformly there: for any sequence of thresholds in that neighborhood, a convergent subsequence satisfies the theorem. In the OU case use \eqref{eq:ou-relative}. The continuous-exponent result (or the same relative approximation at the endpoints) brackets the exact quantile in this neighborhood, since the rates at the endpoints lie on opposite sides of $r$. Thus at $b_d^*$,
\[
 -dr=-dI_{H_d}(b_d^*)-\tfrac12\log d+\log A(H_d,b_d^*)+o(1),
 \quad A(H,b)=\frac{D}{\theta_b\sqrt{2\pi V}}.
\]
In the OU case the $o(1)$ is $O(d^{-1})$. All the quantities in $A$ are bounded above and away from zero there. The rate derivative is uniformly separated from zero, so $|b_d^*-b_{r,d}|=O(\log d/d)$.

On this separated multiplier interval $I''=1/V$ is uniformly bounded, and $A$ is locally Lipschitz in $b$, uniformly in the proportions. Indeed the log-potential response in Proposition~\ref{prop:global-response} bounds the derivative of $g$; the covariance derivative then bounds the kernel and diagonal changes of $K$. For the determinant, the trace series and the Hilbert--Schmidt inequality bound the change of each $j\ge2$ trace by $Cj q_0^{j-2}$ times the kernel change, and the first trace by the diagonal change. The inverse in $V$ is bounded by the gap. These give the Lipschitz bound for $\log A$, since $\theta,V$ stay positive. Taylor expansion now gives
\[
 d\theta_{r,d}(b_d^*-b_{r,d})=\tfrac12\log d-\log A_{r,d}+O((\log d)^2/d)+o(1),
\]
with the final $o(1)$ replaced by $O(d^{-1})$ for OU data, proving the two remainders.
\end{proof}

At any finite dimension for which a relative bound $|P_d/P_{0,d}-1|\le e_d<1$ is available throughout a quantile bracket, the exact quantile also satisfies
\[
 \frac{e^{-dr}}{1+e_d}\le P_{0,d}(b_d^*)\le\frac{e^{-dr}}{1-e_d}.
\]
On a separated saddle interval the Lipschitz bound for $\log A$ and $I'_H=-\theta_b$ imply that $\log P_{0,d}$ is strictly increasing for sufficiently large $d$. The preceding displays therefore locate the quantile, including the essential $\log d/(2d\theta_b)$ displacement from the rate-only threshold.

\section{Temporal covariance and power comparisons}
\label{sec:ou-comparisons}

Changing temporal correlation while keeping the marginal variances fixed generally does not give a positive-semidefinite order of covariance operators: a nonzero positive continuous kernel cannot have identically zero diagonal. The comparison below instead uses the determinants of finite sets of observations. It gives both an exponential-scale and a finite-sample result.

\begin{theorem}[Attenuating temporal covariances]
\label{thm:attenuation}
For each type replace $C_{k,0}$ by
\begin{equation}
 C'_{k,0}(s,t)=C_{k,0}(s,t)R_k(s,t),
 \label{eq:attenuation-kernel}
\end{equation}
where $R_k$ is a continuous positive-semidefinite kernel with $R_k(t,t)=1$. Assume the resulting prior laws satisfy the continuous-path assumptions of Section~\ref{sec:global-structure}. Keep powers and proportions fixed. The two models have the same typical service $b_0$.
\begin{enumerate}
\item For every common internal threshold $0<b<b_0$, $I_{H'}(b)\ge I_H(b)$.
\item At any finite path count and every $\lambda>0$,
\begin{equation}
 \E e^{-\lambda S'_d}\le\E e^{-\lambda S_d},\qquad
 \E S'_d=\E S_d,\qquad \operatorname{Var}(S'_d)\le\operatorname{Var}(S_d).
 \label{eq:laplace-comparison}
\end{equation}
\end{enumerate}
If, for some type of positive proportion (respectively, positive count in part 2), there are distinct $s,t$ with $C_{k,0}(s,t)\ne0$ and $|R_k(s,t)|<1$, the rate, Laplace-transform and variance inequalities are strict. Thus the two finite-sample CDFs cross: there are positive thresholds with each of the two strict orderings $P'_d(b)<P_d(b)$ and $P'_d(b)>P_d(b)$. If the limiting proportions in the large-deviation theorem are positive, strict attenuation gives for every fixed internal $b$
\begin{equation}
 \lim_{d\to\infty}\frac1d\log\frac{P'_d(b)}{P_d(b)}
 =-(I_{H'}(b)-I_H(b))<0.
 \label{eq:strict-ratio}
\end{equation}
Moreover, suppose the counts have convergent proportions and the strictness condition holds for a type of positive limiting proportion. Let $b_{0,d}$ denote the common typical service at the actual proportions. For each $c>0$, for all sufficiently large $d$ there is a crossing of the two CDFs in
$(b_{0,d}-c/\sqrt d,b_{0,d}+c/\sqrt d)$.
Under the same strict attenuation condition, for fixed rational proportions let $\sigma,\mu,\gamma$ and $\sigma',\mu',\gamma'$ be the coefficients of Proposition~\ref{prop:edgeworth} for the two models. Then $\mu'=\mu$, $\sigma'<\sigma$, and there exists a crossing $\widehat b_d$ with
\begin{equation}
 \widehat b_d=b_0+\frac{y_*}{d}+o(d^{-1}),\qquad
 y_* = \mu-\frac{\gamma'/(6\sigma'^3)-\gamma/(6\sigma^3)}{1/\sigma'-1/\sigma}.
 \label{eq:crossing-displacement}
\end{equation}
\end{theorem}

\begin{proof}
\emph{Determinant comparison.} For positive-definite covariance matrix $C$ and correlation matrix $R$, the Schur (entrywise) product satisfies
\begin{equation}
 \det(C\circ R)\ge\det C.
 \label{eq:schur-det}
\end{equation}
This is a consequence of the classical spectral majorization for Schur products \cite{BapatSunder1985}. A direct proof useful here is as follows. Represent $C$ and $R$ by Gram vectors $u_j$ and $v_j$, with $\|v_j\|=1$. The Schur product is the Gram matrix of $u_j\otimes v_j$. The distance of $u_j\otimes v_j$ to the span of the preceding tensor vectors is at least its distance to $\operatorname{span}\{u_i:i<j\}\otimes\operatorname{span}\{v_i:i\le n\}$, which is exactly the distance of $u_j$ to the preceding $u_i$. Taking the product of the squared distances gives \eqref{eq:schur-det}. It holds for every principal submatrix as well. For a two-by-two submatrix indexed by times $s,t$, the determinant increases by $C(s,t)^2(1-R(s,t)^2)$.

\emph{Large-deviation costs.} For a fixed continuous $g>0$, the Gaussian quadratic transform is $Z_k(g)=\Delta_k(g)^{-1/2}$, where
\[
 \Delta_k(g)=\det_F(I+2p_kM_{\sqrt g}C_{k,0}M_{\sqrt g}).
\]
The Fredholm expansion is
\[
 \Delta_k(g)=1+\sum_{j\ge1}\frac{(2p_k)^j}{j!}
 \int_{[0,T]^j}\prod_{\ell=1}^jg(t_\ell)
     \det[C_{k,0}(t_i,t_\ell)]_{i,\ell=1}^j\,dt_1\cdots dt_j.
\]
All terms are nonnegative, and the series is dominated by the exponential series from the diagonal bounds (Hadamard's inequality). By \eqref{eq:schur-det}, $\Delta'_k(g)\ge\Delta_k(g)$, strictly for the indicated type, since the strict two-time determinant difference persists on a set of positive Lebesgue measure by continuity. Hence $-\log Z'_k(g)\ge-\log Z_k(g)$, with strict inequality for that type for every $g>0$. In the variational expression for $\Lambda$ preceding Proposition~\ref{prop:global-response}, the only term that changes is $-\sum_k\pi_k\log Z_k(g)$. Evaluate the new objective at its own minimizer to conclude $\Lambda_{H'}(\theta)\ge\Lambda_H(\theta)$, strictly under the stated condition, for every $\theta>0$. Finally,
\[
 I_{H'}(b)=\sup_{\theta\ge0}(\Lambda_{H'}(\theta)-\theta b)
 \ge\Lambda_{H'}(\theta_b)-\theta_b b\ge I_H(b).
\]
Here $\theta_b>0$ is the original saddle. Strictness follows from the preceding strict inequality. Equation \eqref{eq:strict-ratio} now follows from Theorem~\ref{thm:continuous-exponent}.

\emph{Finite-sample Laplace transforms.} Use positive time quadrature with nodes $t_j$ and weights $a_j$. For independent gamma variables $G_j$ of shape $\lambda a_j$ and rate one, the scalar gamma transform gives
\[
 \exp\left\{-\lambda\sum_j a_j\log(1+\rho q_d(t_j))\right\}
 =\E_G\exp\left\{-\rho\sum_jG_j q_d(t_j)\right\}.
\]
With $D_G=\operatorname{diag}(G_j)$ and $C_{k,N}$ the prior covariance matrix at the nodes, the expectation over the paths, conditional on $G$, is
\[
 \prod_k\Delta_{k,N}(G)^{-n_k/2},\qquad
 \Delta_{k,N}(G)=\det\left(I+\frac{2\rho p_k}{d}D_G^{1/2}C_{k,N}D_G^{1/2}\right).
\]
Expand each determinant into its principal minors. Each weighted minor increases under \eqref{eq:attenuation-kernel} by \eqref{eq:schur-det}. Since all $G_j\ge0$, the product of inverse determinant powers decreases. Taking expectations over the gamma variables gives the ordering for the grid services. Continuous paths and the bound $0\le e^{-\lambda S_{d,N}}\le1$ permit passage to the continuous integral by dominated convergence.

Strictness also survives the grid limit. If the condition of the theorem holds for type $k$, choose disjoint relative open time intervals $U,V$ of positive length and $\delta>0$ such that $C_{k,0}(s,t)^2(1-R_k(s,t)^2)\ge\delta$ on $U\times V$. The contribution of the two-by-two principal minors to $\Delta'_{k,N}-\Delta_{k,N}$ is then at least
\[
 (2\rho p_k/d)^2\delta\left(\sum_{t_j\in U}G_j\right)\left(\sum_{t_j\in V}G_j\right).
\]
For fine uniform midpoint grids, the two sums are independent gamma variables with shapes tending to $\lambda|U|$ and $\lambda|V|$, respectively. The sum over remaining nodes is independent and its shape is bounded by $\lambda T$. With a probability bounded below uniformly on these grids, both displayed sums lie in $[1,2]$ and the remaining sum is at most one. (Choose intervals with nonempty complement; the limiting gamma shapes are then positive.) On this event all determinants lie in $[1,B]$ for some finite $B$ independent of the grid, by $\det(I+A)\le e^{\operatorname{tr}A}$ for $A\succeq0$ and $\sum_j G_j\le5$. The indicated determinant increases by at least a fixed positive amount. For $n_k>0$, $x\mapsto x^{-n_k/2}$ has a strictly negative derivative bounded away from zero in absolute value on $[1,B]$. Multiplying by the other inverse determinant powers, which are bounded below by positive constants, gives a strictly positive grid-independent difference on an event of grid-independent positive probability. This proves the strict continuous Laplace ordering.

\emph{Means, variances and crossing.} The same diagonal covariance gives the same one-time laws of $q_d(t)$ and hence the same service mean. For $a_k=\rho p_k/d$ the two-time transform is
\begin{align*}
 \E e^{-u\rho q_d(s)-v\rho q_d(t)}
 &=\prod_k\big[(1+2u a_kC_{k,0}(s,s))(1+2v a_kC_{k,0}(t,t))\\
 &\hspace{5em}-4uv a_k^2 C_{k,0}(s,t)^2\big]^{-n_k/2},\qquad u,v>0.
\end{align*}
It decreases under the attenuation, strictly on $U\times V$ under the strictness condition. The product of the two marginal transforms is unchanged. Using $\log(1+x)=\int_0^\infty e^{-u}(1-e^{-ux})du/u$, the covariance of the two log-services is the integral of $e^{-u-v}/(uv)$ times the difference between the joint and product transforms. This identity follows by truncating the integrals and using convergence in $L^2$, since the logarithm is bounded by $x$ for $x\ge0$. The difference of the transforms is nonnegative by the product formula. The covariance therefore decreases; integrating over $(s,t)\in[0,T]^2$ gives the variance inequality, strictly under the condition. With equal finite means and different variances, the two nonnegative service laws cannot have CDFs ordered at every positive threshold: right continuity would make the integral of the survival functions strictly ordered, contradicting equality of the means. This proves the crossing assertion.

To locate a crossing, linearize the service about $m_{0,d}=\sum_k\pi_{k,d}p_k C_{k,0}(t,t)$. The Taylor remainder after multiplying by $\sqrt d$ has mean absolute value at most $C\sqrt d\int\E(q_d-m_{0,d})^2dt=O(d^{-1/2})$, by the bounded second derivative of the logarithm and Wick's identity. The linear term is a sum of independent centered Gaussian quadratic forms divided by $\sqrt d$, with uniformly bounded unscaled third absolute moments. Its variance tends to
\[
 \sigma^2=2\sum_k\pi_kp_k^2\iint_{[0,T]^2}
 w_0(s)w_0(t)C_{k,0}(s,t)^2ds\,dt,\qquad
 w_0=\frac{\rho}{1+\rho m_0}.
\]
By the triangular-array central limit theorem
$\sqrt d(S_d-b_{0,d})\Rightarrow N(0,\sigma^2)$. The same argument applies to $S'_d$ with variance $\sigma'^2$. Because $w_0$ is unchanged, $R_k^2\le1$, and the strict attenuation persists on a set of positive measure for a type with $\pi_k>0$, we have $0<\sigma'^2<\sigma^2$. Thus for fixed $x\ne0$,
\[
 P'_d(b_{0,d}+x/\sqrt d)-P_d(b_{0,d}+x/\sqrt d)
 \longrightarrow \Phi(x/\sigma')-\Phi(x/\sigma).
\]
The right side is negative for $x=-c$ and positive for $x=c$. Both CDFs are continuous by the argument in Section~\ref{sec:quantiles}. The intermediate value theorem gives a crossing between these two thresholds. For fixed rational proportions, the common diagonal kernels give the same $\mu$ in Proposition~\ref{prop:edgeworth}. Apply its uniform expansion at $b=b_0+y/d$. Taylor expansion of the two normal CDFs gives, uniformly for bounded $y$,
\[
 \begin{aligned}
 \sqrt d\{P'_d(b_0+y/d)-P_d(b_0+y/d)\}
  =\phi(0)\bigg[&(y-\mu)(1/\sigma'-1/\sigma)\\
  &+\frac{\gamma'}{6\sigma'^3}-\frac{\gamma}{6\sigma^3}\bigg]+o(1).
 \end{aligned}
\]
The limiting affine function has positive slope and unique zero $y_*$. Uniform convergence on a fixed compact interval around $y_*$ and continuity give a zero for each large $d$; every such zero in this interval converges to $y_*$. This proves \eqref{eq:crossing-displacement} without requiring uniqueness of the finite-sample crossing.
\end{proof}

These comparisons use an order on the entire covariance kernel, not just a comparison of correlation values at a single time lag. In particular, the hypothesis requires that the multiplier $R_k$ itself be a correlation kernel. Products with $R(s,t)=e^{-a|s-t|}$ or $R(s,t)=e^{-a(s-t)^2}$, $a>0$, give distinct examples, the latter generally destroying the Markov property of an OU process.

\subsection{Sensitivity to an OU correlation time}

For $C_{k,0}(s,t)=e^{-\kappa_k|s-t|}$ with $\kappa_k>0$, increasing $\kappa_k$ is an instance of \eqref{eq:attenuation-kernel}. The rate also admits an explicit derivative in this case.

\begin{corollary}[OU rate sensitivity]
\label{thm:ou-clock}
Fix $T,\rho,\pi,p$ and a strictly internal $b$. The rate is differentiable in the OU decay rates, and
\begin{equation}
 \frac{\partial I_H(b)}{\partial\kappa_k}=-\pi_k\partial_{\kappa_k}\log Z_k(g_{\theta_b})>0,
 \label{eq:clock-envelope}
\end{equation}
where differentiation of $Z_k$ keeps the potential fixed. For a reliability level $r>0$ admitting an internal root $I_H(b_r)=r$, $\partial_{\kappa_k}b_r=(\partial_{\kappa_k}I_H(b_r))/\theta_{b_r}>0$.
\end{corollary}

\begin{proof}
For ordered distinct times the OU Markov covariance has determinant $\prod_{j=1}^{n-1}(1-e^{-2\kappa(t_{j+1}-t_j)})$, by successive Gaussian innovations. For $n\ge2$ its derivative in $\kappa$ is strictly positive and at most $2\sum_j(t_{j+1}-t_j)\le2T$. The Fredholm expansion used above thus differentiates term by term, by domination by an exponential series. Its two-time term ensures $-\partial_\kappa\log Z_k(g)>0$ for fixed $g>0$. Differentiability of the tilted covariance in $\kappa$ follows by differentiating its continuous-kernel resolvent equation; the inverse in that equation is bounded. The fixed-point derivative is invertible by the contraction, and $B'(\theta)<0$. The implicit-function theorem therefore gives differentiability of the saddle. Stationarity in $g$ and $\Lambda'(\theta_b)=b$ cancel the saddle derivative terms in $I_H(b)=\Lambda(\theta_b)-\theta_b b$, giving \eqref{eq:clock-envelope}. The last formula follows from $\partial_b I=-\theta_b$.
\end{proof}

The exponential-scale advantage holds at every fixed internal threshold as $d$ grows; it does not contradict the crossing of the CDFs at each finite $d$.

If both models satisfy the compact-family assumptions of Theorem~\ref{thm:ou-relative}, that result also supplies a sufficient finite-count comparison: with a common relative error bound $e_d<1$, $P'_d(b)<P_d(b)$ whenever $(1+e_d)P'_{0,d}(b)<(1-e_d)P_{0,d}(b)$. At a fixed internal threshold this condition eventually holds under strict attenuation, because the rate gap in \eqref{eq:strict-ratio} dominates the bounded prefactor ratio.

\subsection{Power dispersion for a common temporal kernel}
Here the powers vary while all types share the same temporal kernel.
The comparison in Theorem~\ref{thm:attenuation} instead varies temporal dependence at fixed powers.

\begin{proposition}[Power dispersion with a common temporal kernel]
\label{prop:power-dispersion}
Suppose all types have the same nonzero continuous covariance kernel
$C_0$. Let $H_1,H_2$ be two finite probability laws of positive powers,
with the same mean, such that $H_1\preceq_{\rm cx}H_2$ in convex order.
Then at every common strictly internal threshold,
\[
 I_{H_2}(b)\le I_{H_1}(b).
\]
The inequality is strict if the two finite laws differ. In particular,
equal powers maximize this variational reliability exponent among
allocations with the same mean and a common temporal covariance.
For any common internal reliability level, $b_r(H_2)\le b_r(H_1)$.
\end{proposition}

\begin{proof}
For fixed $g>0$, set $A=C_0^{1/2}M_gC_0^{1/2}$ and
\[
 f_g(p)=\tfrac12\log\det_F(I+2pA).
\]
It is strictly concave, because
$f_g''(p)=-2\operatorname{tr}(A^2(I+2pA)^{-2})<0$.
Thus $E_{H_2}f_g\le E_{H_1}f_g$. To see strictness for distinct finite laws, take a compact interval $[a,c]$ containing their supports and use $f_g(p)=f_g(a)+f_g'(a)(p-a)+\int_a^c f_g''(v)(p-v)_+dv$. Convex order makes the difference of the expectations of $(p-v)_+$ nonnegative; for distinct laws it is positive on a nonempty interval (these continuous functions determine the laws). The means agree and $f_g''<0$, so the inequality is strict. In the variational formula for $\Lambda_H(\theta)$,
the only dependence on $H$ is $E_Hf_g$. Taking the infimum over $g$
gives $\Lambda_{H_2}(\theta)\le\Lambda_{H_1}(\theta)$.
Strictness follows by evaluating the $H_2$ objective at the minimizer
for $H_1$; the fixed-point potential is strictly positive.
For the $H_2$ saddle $\theta_2>0$,
\[
 I_{H_2}(b)=\Lambda_{H_2}(\theta_2)-\theta_2b
 <\Lambda_{H_1}(\theta_2)-\theta_2b\le I_{H_1}(b),
\]
when the laws differ. The last inequality uses the exact variational
representation $I_H(b)=\sup_{\theta\ge0}(\Lambda_H(\theta)-\theta b)$.
Monotonicity in $b$ gives the threshold conclusion.
\end{proof}

\section{Evaluation and illustrations}
\label{sec:illustration}

The saddle and Fredholm factor can be evaluated by Nystr\"om quadrature \cite{Bornemann2010}. On an $N$-cell midpoint grid of width $h=T/N$, write $G=h\operatorname{diag}(g_j)$ and factor the prior matrix $C_{k,N}=L_kL_k^{\mathsf T}$. One fixed-point iteration uses
\[
 \begin{aligned}
 S_{k,N}&=L_k(I+2p_kL_k^{\mathsf T}GL_k)^{-1}L_k^{\mathsf T},\\
 g_j^{\rm new}&=\frac{\theta\rho}{1+\rho\sum_k\pi_kp_k(S_{k,N})_{jj}},
 \qquad -\log Z_{k,N}=\tfrac12\log\det(I+2p_kL_k^{\mathsf T}GL_k).
 \end{aligned}
\]
No prior inverse is needed. The contraction and \eqref{eq:global-residual} hold independently of $N,d$. Solve $B_N(\theta)=b$ by bracketing. At its solution $I_N=-\sum_k\pi_k\log Z_{k,N}-h\sum_jg_jm_{g,j}$. To avoid cancellation for small $\theta$, compute it as $\tfrac12\sum_{k,j}\pi_k[\log(1+\lambda_{kj})-\lambda_{kj}/(1+\lambda_{kj})]$, using its expansion near zero, where $\lambda_{kj}$ are the eigenvalues of $2p_kL_k^{\mathsf T}GL_k$. Proposition~\ref{prop:edgeworth} similarly avoids cancellation in $D/u-1/r$. Form
\[
 \begin{gathered}
 (K_N)_{ij}=2\theta h\,w_iw_j\sum_k\pi_kp_k^2(S_{k,N})_{ij}^2,\\
 D_N=\det(I-K_N)^{-1/2},\qquad
 V_N=\tfrac h\theta\mathbf1^{\mathsf T}K_N(I-K_N)^{-1}\mathbf1.
 \end{gathered}
\]
The covariance iteration costs $O(mN^3)$ and the final curvature calculation $O(N^3)$ operations. Section~\ref{sec:continuous-exponent} proves convergence of $I_N$; the uniform kernel convergence and the trace and resolvent series with the common gap give convergence of $D_N,V_N$, as in Theorem~\ref{thm:sharp}.

\subsection{Continuous Riccati evaluation for OU types}
\label{sec:riccati}

For OU types, the scalar Riccati and variance equations of Section~\ref{sec:ou-quantitative} directly give $C_{k,g}(t,t)=c_k(t)$ and $Z_k(g)=e^{A_k(0)}/\sqrt{1+r_k(0)}$ for the fixed-point iteration. The following system computes the prefactor with state size depending only on $m$.

\begin{proposition}[Riccati representation of the prefactor]
\label{prop:riccati}
Use the OU radial Gaussian reference and functions $f_k,\tau_k$ in Section~\ref{sec:ou-quantitative}, and write $f=(f_1,\ldots,f_m)^{\mathsf T}$, $A=2\operatorname{diag}(\tau_k')$, $Q=2\operatorname{diag}(\tau_k')$. (Thus the reference vector satisfies $dZ=-AZdt+Q^{1/2}dW$.) Define the symmetric matrix $R(t)$ and scalar $c(t)$ by
\begin{equation}
 R'=AR+RA-RQR-\theta ff^{\mathsf T},\qquad c'=-\tfrac12\operatorname{tr}(QR),
 \qquad R(T)=0,\quad c(T)=0.
 \label{eq:curvature-riccati}
\end{equation}
These equations have a bounded solution on $[0,T]$, with $R(t)\succeq0$ and $I-R(0)/2\succeq(1-q_0)I$. Let $F(t)$ and $\overline Q(t)$ be the augmented matrices
\[
 F=\begin{pmatrix}-A+QR&0\\ f^{\mathsf T}&0\end{pmatrix},\qquad
 \overline Q=\begin{pmatrix}Q&0\\0&0\end{pmatrix}.
\]
Solve $S'=FS+SF^{\mathsf T}+\overline Q$, with $S(0)=\operatorname{diag}((2I-R(0))^{-1},0)$. Then
\begin{equation}
 D=e^{c(0)}\det(I-R(0)/2)^{-1/2},\qquad V=S_{m+1,m+1}(T).
 \label{eq:riccati-prefactor}
\end{equation}
\end{proposition}

\begin{proof}
For the reference $B=f^{\mathsf T}Z$, the weight $\exp(\theta\int B^2/2)$ has mass $D$ by the curvature gap. Given $Z(t)=z$, the future $B$ is Gaussian with linear-transition mean $M_tz$ and covariance $C_t^{\rm cond}$ dominated by the restriction of $K/\theta$. Completing the Gaussian square gives the conditional future weight
\[
 h(t,z)=\det_F(I-\theta C_t^{\rm cond})^{-1/2}
   \exp\big(\tfrac\theta2\langle M_tz,(I-\theta C_t^{\rm cond})^{-1}M_tz\rangle\big).
\]
The gap $\theta\|C_t^{\rm cond}\|\le q_0$ bounds both coefficients of $h(t,z)=e^{c(t)+z^{\mathsf T}R(t)z/2}$, namely $R(t)=\theta M_t^*(I-\theta C_t^{\rm cond})^{-1}M_t$ and $c(t)$, and both are nonnegative.

For the differential equations, use a left-endpoint quadrature of the future integral. Backward Gaussian conditioning over a step $h$ with transition matrix $\Phi=I-Ah+o(h)$ and covariance $\Xi=Qh+o(h)$ gives
\[
 R_-=h\theta ff^{\mathsf T}+\Phi^{\mathsf T}R_+(I-\Xi R_+)^{-1}\Phi,\qquad
 c_-=c_+-\tfrac12\log\det(I-\Xi R_+).
\]
The sampled unconditional covariance of $B$ has a strict norm gap on fine grids by uniform convergence of the weighted rows of $K$. Conditional covariance is smaller, so the transform formula bounds the grid coefficients and a slightly higher exponential moment for each fixed initial mean. Uniform integrability and path continuity give convergence of the conditional integrals. The recurrence has $O(h)$ increments uniformly and first-order consistency with \eqref{eq:curvature-riccati}; its piecewise linear interpolants are equicontinuous. Passing to the integral equations proves bounded solvability and local Lipschitz continuity of the vector field gives uniqueness.

Since $Z(0)\sim N(0,I/2)$, the unconditional covariance is $C_0^{\rm cond}+M_0M_0^*/2$. Conjugating its gap by $(I-\theta C_0^{\rm cond})^{-1/2}$ gives
\[
 I-\tfrac{\theta}{2}(I-\theta C_0^{\rm cond})^{-1/2}M_0M_0^*
 (I-\theta C_0^{\rm cond})^{-1/2}\succeq(1-q_0)I,
\]
using $(I-\theta C_0^{\rm cond})^{-1}\succeq I$. The subtracted operator has the same nonzero eigenvalues as $R(0)/2$, giving the initial gap. Integrating $h(0,z)$ against $N(0,I/2)$ gives $D$.

The normalized product of the accrued exponential and $h(t,Z_t)$ is the density martingale by conditional expectation. It\^o's formula and localization give drift correction $Q\nabla\log h=QRZ$, unchanged diffusion and initial covariance $(2I-R(0))^{-1}$. Augmentation by $\int_0^t f^{\mathsf T}Z$ gives the stated Lyapunov equation; the terminal integral variance is $V$.
\end{proof}

The scalar and matrix ODEs have respectively $O(m)$ and $O(m^2)$ unknowns. On bounded sets containing exact and approximate solutions, an endpoint error $\delta$ and differential residual $\epsilon$ give error at most $e^{LT}(\delta+T\epsilon)$ by Gronwall, with vector-field Lipschitz bound $L$. The initial matrix gap controls determinant and inverse evaluation; inexact variances enter the fixed point via \eqref{eq:global-residual}.

\subsection{Accuracy of the continuous saddle computation}
\label{sec:quadrature-accuracy}

Since the rate is multiplied by $d$, we need quantitative time-discretization accuracy.

\begin{proposition}[Midpoint saddle and prefactor error]
\label{prop:midpoint-error}
For fixed Gaussian types and powers, let a threshold $b=B(\theta)$ have its saddle in a compact positive multiplier interval. Use actual proportions. Let $h=T/N$ and set
\[
 \omega(h)=\max_k\sup_{|t-t'|+|s-s'|\le h}
       |C_{k,0}(t,s)-C_{k,0}(t',s')|.
\]
The midpoint-grid saddle at the same threshold satisfies, for all fine grids,
\begin{equation}
 |\theta_N-\theta|+|I_N-I|+|\log D_N-\log D|+|V_N-V|
 \le C\,\omega(h).
 \label{eq:midpoint-general}
\end{equation}
For OU types on a compact set of positive powers and rates, the right-hand side can be replaced by $Ch^2$, uniformly on compact positive multiplier intervals. All constants may be chosen uniformly in the proportions for the given finite type list. Consequently, the logarithmic error in evaluating $P_{0,d}$ is $O((d+1)\omega(h))$, or $O((d+1)h^2)$ for OU kernels, in exact arithmetic.
\end{proposition}

\begin{proof}
At fixed $\theta$ the resolvent satisfies
\[
 C_{k,g}(t,s)=C_{k,0}(t,s)-2p_k\int_0^T C_{k,0}(t,v)g(v)C_{k,g}(v,s)\,dv.
\]
Subtracting in the first argument, then using symmetry and $|C_{k,g}|\le M_k$, bounds the modulus of $C_{k,g}$ by a constant times the prior modulus; the same holds for $g,w$. The midpoint error in this integral is $O(\omega(h))$ uniformly at node pairs. For $G=h\operatorname{diag}(g(t_j))$, the exact nodal covariance matrix $S^C$ thus satisfies $(I+2p_k C_{k,N}G)S^C=C_{k,N}+E_N$, with entrywise $E_N=O(\omega(h))$. The inverse $I-2p_k S_{k,N}G$ has maximum absolute row sum at most $1+2p_kTM_k\|g\|_\infty$, so the grid-tilted covariance at these potentials has entrywise error $O(\omega(h))$. The log-potential residual has the same order. Contraction and the resolvent derivative give that order for the actual grid potential and covariance, and service quadrature then gives $|B_N-B|\le C\omega(h)$.

On compact positive multiplier intervals $V$ dominates the positive variance lower bound in Lemma~\ref{lem:central-half-estimates}, uniformly in proportions. Since $B'=-V$, a small bracket gives $|\theta_N-\theta|\le C\omega(h)$, preserving the previous estimates at the roots. For the normalizer, differentiation in the scalar multiplier of the prescribed $g$ gives
\[
 -\log Z_k(g)=p_k\int_0^1\int_0^T g(t)C_{k,sg}(t,t)\,dt\,ds.
\]
The grid identity is identical with time quadrature. The nodal error argument applies uniformly to $sg$, $0\le s\le1$, which has the same modulus bound. Thus $I=-\sum_k\pi_k\log Z_k-\int gm_g$ has error $O(\omega(h))$.

The curvature kernel has the same modulus and nodal error, including the factor $h$ for the grid matrix. With $u=(I-K)^{-1}1$, $u=1+Ku$ and $\|u\|_\infty\le(1-q_0)^{-1}$ give the same modulus for $u$. Substitution at the nodes gives residual $O(\omega(h))$ and inversion by the grid row bound controls $u$ and $V=\theta^{-1}\int(u-1)$.

For $\log D=\tfrac12\sum_{j\ge1}\operatorname{tr}(K^j)/j$, the sampled continuous matrix $\widehat K_N=(hK(t_i,t_j))_{ij}$ has row bound $q_1<1$ on fine grids for fixed $q_1>q_0$. In the cyclic trace integral for $K^j$, $j\ge2$, replace one integral at a time by midpoint sums. Only two factors depend on the replaced variable (a square for $j=2$), giving error $O(\omega(h))$; the remaining chain of $j-2$ factors integrates or sums to at most $Tq_1^{j-2}$ by successive row bounds. The trace error is at most $Cj q_1^{j-2}\omega(h)$. Also $\|\widehat K_N-K_N\|_{\rm HS}=O(\omega(h))$; telescope the powers under the trace using this difference, one other Hilbert--Schmidt factor and the remaining operator norms for the same bound. The first trace uses just diagonal quadrature. Summing proves the log-determinant estimate.

For OU types, the backward Riccati and forward variance equations give bounded solutions and coefficients on compact data. Since $g$ is continuous, $c_k\in C^1$, hence $g\in C^1$ by the fixed point; differentiation gives $c_k,g\in C^2$ with uniform derivative bounds. The same bounds hold for variances at $sg$, $0\le s\le1$. The tilted covariance $\sqrt{c_k(t)c_k(v)}e^{-|\tau_k(t)-\tau_k(v)|}$ and the prior are continuous, piecewise $C^2$ in either argument, with uniformly bounded one-sided derivatives and just one break at equality. A continuous piecewise $C^2$ function of this kind, with finitely many breaks, has midpoint error $O(h^2)$ (smooth cells contribute $O(h^3)$ each, break cells $O(h^2)$). The resolvent integrand has at most two breaks; service and diagonal normalizer integrands are smooth. Thus their previous errors are $O(h^2)$. The same piecewise structure holds for $K$. Splitting $u=1+Ku$ at $s=t$ shows $u\in C^2$ with bounded derivatives: boundary terms cancel on first differentiation by continuity and the second contains only the bounded jump of the first derivative times $u(t)$. Hence quadrature of $K(t_i,s)u(s)$ and $\int u$ has the improved order. Each product in the trace argument has at most two breaks, so its bound also holds with $h^2$. The logarithm of \eqref{eq:effective-centre}, with $\theta,V$ bounded away from zero, proves the probability error.
\end{proof}

Log-potential and service-root residuals of order $\epsilon$ give additional saddle and prefactor errors $O(\epsilon)$ by contraction and strict response. For OU paths, $h=O(d^{-1})$, $\epsilon=O(d^{-2})$ keep log-probability error on the analytic $O(d^{-1})$ scale; for general kernels require $d\omega(h),d\epsilon\to0$. Bracketing and contraction cost $O(mN^3(1+|\log\epsilon|)^2)$ operations before the curvature calculation. These are exact-arithmetic estimates; certified numerical bounds also require control of roundoff.

\subsection{A continuous-probability bound for reversible OU paths}
\label{sec:continuous-grid-bound}

Exact grid transitions also give finite-count bounds for the continuous event.

\begin{proposition}[Trapezoidal probability transfer]
\label{prop:finite-grid-transfer}
Let the independent prior paths be stationary, variance-one OU processes of rates $\kappa_i$, with weights $p_i$. For a partition $0=t_0<\cdots<t_M=T$ of widths $h_j$, let $S_{d,h}$ be the composite trapezoidal rule applied to $\log(1+\rho q_d(t))$ and let $F_h(v)=\Prb(S_{d,h}\le v)$. Put
\[
 K_p=\max_i\kappa_i p_i,\qquad v_h=\frac{\rho K_p}{6d}\sum_j h_j^3.
\]
For every $e>0$ each one-sided integration error has probability at most $\delta_h(e)=\exp(-e^2/(2v_h))$. In particular
\begin{equation}
 F_h(b-e)-\delta_h(e)\le P_d(b)\le F_h(b+e)+\delta_h(e).
 \label{eq:finite-grid-transfer}
\end{equation}
For a uniform mesh of width $h$, the exponent in $\delta_h(e)$ is $-3de^2/(\rho K_pTh^2)$.
\end{proposition}

\begin{proof}
Put $f(x)=\log(1+\rho d^{-1}\sum_i p_ix_i^2)$, $w_h(t)=t-(t_j+t_{j+1})/2$ on each panel, and $\mathcal L=\sum_i\kappa_i(\partial_{ii}-x_i\partial_i)$. Integration by parts and It\^o's formula give $H=S_{d,h}-S_d=\int w_h\,df(X_t)=A+J$, where $A=\int w_h\,dM$, $J=\int w_h\mathcal Lf(X_t)dt$. The reversed stationary OU process $Y_s=X_{T-s}$ has the same Gaussian law and generator. The reversed integration-by-parts identity gives $-H=\widetilde A+J$, $\widetilde A=\int w_h(T-s)d\widetilde M_s$ in its own filtration; the drifts agree by change of variables. Thus $2H=A-\widetilde A$, the reversible forward--backward decomposition \cite{LyonsZheng1988,AltmeyerChorowski2018}. The quadratic variation rate satisfies
\[
 \Gamma_f(x)=\frac{8\rho^2}{d^2}\frac{\sum_i\kappa_i p_i^2x_i^2}{(1+\rho q(x))^2}
 \le\frac{8\rho K_p}{d}\frac{\rho q(x)}{(1+\rho q(x))^2}\le\frac{2\rho K_p}{d}.
\]
Since $\int w_h^2=\sum_jh_j^3/12$, each martingale has quadratic variation at most $v_h$ and exponential moment at most $e^{u^2v_h/2}$ for real $u$. By Cauchy--Schwarz, $\E e^{uH}\le(\E e^{uA}\E e^{-u\widetilde A})^{1/2}\le e^{u^2v_h/2}$; optimizing Chernoff gives the two tails and sandwiching gives \eqref{eq:finite-grid-transfer}.
\end{proof}

For completeness, at fixed $d$ the trapezoidal error has the exact leading mean-square coefficient on uniform grids:
\begin{equation}
 \lim_{h\downarrow0} h^{-2}\E(S_d-S_{d,h})^2=\frac{T}{12}\E\Gamma_f(X_0).
 \label{eq:trapezoid-ms}
\end{equation}
Indeed $\mathcal Lf$ is bounded and continuous for fixed parameters (its quadratic numerators are controlled by the denominators), so $\mathcal Lf(X_t)$ is mean-square continuous. In $J$, subtract $\mathcal Lf(X_{t_j})$ on each panel using the zero integral of $w_h$. Minkowski and $\int|w_h|=Th/4$ give $\|J\|_2=o(h)$. Isometry gives $\E A^2=(Th^2/12)\E\Gamma_f(X_0)$, proving \eqref{eq:trapezoid-ms}; see the general limit results of \cite{Altmeyer2023}.

For a type of $n$ identical OU coordinates, the sum of squares $R_j$ at the nodes has the exact transition
\begin{equation}
 R_{j+1}\ \stackrel{\rm law}{=}\ (1-a_j^2)\chi^2_n\left(\frac{a_j^2R_j}{1-a_j^2}\right),
 \quad a_j=e^{-\kappa h_j},\quad R_0\sim\chi_n^2 .
 \label{eq:exact-ou-transition}
\end{equation}
by the Gaussian OU transition and rotational invariance. Only one radial variable per type is needed, and \eqref{eq:finite-grid-transfer} accounts for the unsampled intervals.

The frozen-channel endpoint $X_i(t)=Z_i$ (independent standard normals, power $p$) checks the normalization outside the distinct-time positive-definiteness assumption. For $v=(e^{b/T}-1)/(\rho p)\in(0,1)$, the event is $\chi_d^2\le dv$, with rate $(v-1-\log v)/2$ and prefactor $[(1-v)\sqrt{\pi d}]^{-1}$ by the scalar density and endpoint Laplace integral. The tilt has $2Tpg=v^{-1}-1$, $g=\theta w$; the sole curvature eigenvalue is $\lambda=2\theta T w^2p^2v^2$. Hence $D=(1-\lambda)^{-1/2}$, $V=T\lambda/(\theta(1-\lambda))$, and $D/(\theta\sqrt{2\pi V})=1/[(1-v)\sqrt\pi]$, as required. The quantitative OU theorem assumes positive rates.

\subsection{Non-Markov kernels and conditional fluctuations}

Take $T=\rho=1$, equal proportions, powers $(0.6,1.4)$ and covariances
\begin{equation}
 C_{k,0}^{(a)}(t,s)=\exp(-\kappa_k|t-s|-a(t-s)^2),\quad
 (\kappa_1,\kappa_2)=(1,4),\quad a\ge0.
 \label{eq:example-kernel}
\end{equation}
The Gaussian correlation factor preserves positive definiteness by the Schur product; increment variance $O(|t-s|)$ and Gaussian moment bounds give continuity by Kolmogorov's criterion. For $a>0$ these are non-Markov kernels. All have typical service $\log2$.
\begin{table}[t]
\centering
\caption{Midpoint approximations to the rate at $b=0.63$ in \eqref{eq:example-kernel}. The strict ordering persists under grid refinement.}
\label{tab:nonmarkov-rate}
\begin{tabular}{rrrr}
\toprule
$N$ & $a=0$ & $a=1$ & $a=4$\\
\midrule
32  & 0.01327853 & 0.01433758 & 0.01653258\\
64  & 0.01333188 & 0.01439966 & 0.01661470\\
128 & 0.01334531 & 0.01441529 & 0.01663538\\
256 & 0.01334867 & 0.01441920 & 0.01664056\\
\bottomrule
\end{tabular}
\end{table}

Table~\ref{tab:nonmarkov-rate} (log-potential residual below $2\times10^{-13}$) illustrates Theorem~\ref{thm:attenuation}: at the finest grid, changing $a$ from 0 to 4 increases the rate by about $0.00329$.

For $a=1$ and $N=24$ midpoint nodes, fix $\theta=0.2$, giving $b_N\simeq0.66466007$, $I_N\simeq0.00281308$. Here the event is the exact grid event, and the second-moment condition is $2\theta\rho T\max_k p_kM_k=0.56<1$. The conditional law provides a target for importance sampling \cite{ChatterjeeDiaconis2018}, but rate-optimal tilts need not be efficient \cite{GlassermanWang1997}; see mixtures of twists for nonlinear exceedances \cite{ChanLai2007} and dynamic changes of measure in queues \cite{DupuisSezerWang2007}.

To sample the tilted powers, use the Bartlett representation: for $n_k\ge N$ the unnormalized second moment has law $L_{k,g}A_kA_k^{\mathsf T}L_{k,g}^{\mathsf T}$, where $L_{k,g}$ is the Cholesky factor of $S_{k,N}$, $A_k$ is lower triangular with independent standard normals below the diagonal and independent positive diagonal entries whose squares are $\chi^2_{n_k-j+1}$ at row $j$. This follows from successive orthogonal projection of the rows of a standard Gaussian $N\times n_k$ matrix. The sample likelihood is $\exp\{d(h\sum_jg_jq_d(t_j)+\sum_k\pi_k\log Z_{k,N})\}$. Rescale event weights and probability approximations by $\sqrt d\,e^{dI_N}$. With rescaled weights $W_i$ (zero for misses), deficits $H_i=\theta d(b_N-S_{d,N}^{(i)})$, the conditional mean estimator $\widehat\mu=\sum_iW_iH_i/\sum_iW_i$ has delta-method sampling SE
\[
 \frac{(\sum_i W_i^2(H_i-\widehat\mu)^2)^{1/2}}{\sum_i W_i}.
\]
For conditional variance use the corresponding centered moment; for correlation use the multivariate delta method on its five means.

\begin{table}[t]
\centering
\caption{Exact-grid probability experiment. All probability columns are rescaled by $\sqrt d\,e^{dI_N}$. The sampling standard errors (SE) are for the rescaled probability estimate and the conditional mean of $\theta d(b_N-S_{d,N})$ (in parentheses). Each row uses 400,000 independent grid samples under the tilted law.}
\label{tab:grid-sharp}
\footnotesize
\setlength{\tabcolsep}{4pt}
\begin{tabular}{rrrrrrr}
\toprule
$d$ & Estimate & SE & $P_0$ & $P_{\rm M}$ & $P_{\rm LR}$ & Mean deficit\\
\midrule
128  & 3.5772 & 0.0062 & 5.6916 & 3.4140 & 3.5847 & 0.4803 (.0008)\\
512  & 4.6733 & 0.0102 & 5.6916 & 4.5649 & 4.6553 & 0.7008 (.0013)\\
2048 & 5.3352 & 0.0162 & 5.6916 & 5.2817 & 5.3165 & 0.8727 (.0021)\\
8192 & 5.5659 & 0.0243 & 5.6916 & 5.5726 & 5.5830 & 0.9618 (.0031)\\
\bottomrule
\end{tabular}
\end{table}

Table~\ref{tab:grid-sharp} shows that the leading factor overestimates the probability at smaller counts, while the normal-tail approximations are closer. At $d=8192$ the estimated grid probability is about $6.0\times10^{-12}$ and the deficit mean approaches one. At $d=128$, 400,000 independent prior draws give rescaled estimate $3.5765$ (SE $0.0106$), agreeing with the tilt. The relative second moment divided by $\sqrt d$ is $0.1934,0.1285,0.1032,0.0951$ at the four counts, approaching $0.0916$ from \eqref{eq:is-relative-second-moment}.

For $f(t)=\cos(2\pi t/T)$, set $A_d=h\sum_j f(t_j)w_j\sqrt d\{q_d(t_j)-m_g(t_j)\}$. The fixed-grid version of Theorem~\ref{thm:sharp} predicts conditional variance $\langle f,\Sigma_{N,\perp}f\rangle_h=0.0352741$ with $\langle u,v\rangle_h=h\sum_j u_jv_j$, and zero limiting correlation with $H=\theta d(b_N-S_{d,N})$. Table~\ref{tab:grid-joint} uses the same samples; the variance approaches its limit, and correlation decreases but remains larger than its SE at the largest count.

\begin{table}[b]
\centering
\caption{Conditional projection variance and its correlation with the scaled deficit in the 24-node experiment of Table~\ref{tab:grid-sharp}. Sampling standard errors are in parentheses.}
\label{tab:grid-joint}
\small
\begin{tabular}{rrr}
\toprule
$d$ & $\operatorname{Var}(A_d\mid S_{d,N}\le b_N)$ & $\operatorname{Corr}(A_d,H\mid S_{d,N}\le b_N)$\\
\midrule
128  & 0.03340 (.00011) & 0.0693 (.0020)\\
512  & 0.03461 (.00014) & 0.0549 (.0021)\\
2048 & 0.03505 (.00017) & 0.0430 (.0026)\\
8192 & 0.03517 (.00024) & 0.0218 (.0034)\\
Limit & 0.03527 & 0\\
\bottomrule
\end{tabular}
\end{table}

\subsection{Continuous OU probabilities}

Use \eqref{eq:example-kernel} with $a=0$, equal counts and $b=0.66$. For $d=32,128,512$, draw 50,000 independent prior radial paths by \eqref{eq:exact-ou-transition} on 2048 steps; integrate by the trapezoidal rule. Here $K_p=5.6$, so for the shifts in Table~\ref{tab:ou-continuous}, $\delta_h(e)<10^{-4}$. With $k_-,k_+$ the hit counts below $b-e,b+e$ in $n=50000$ draws, let $L(k),U(k)$ be the $\alpha$ and $1-\alpha$ quantiles of $\operatorname{Beta}(k,n-k+1)$ and $\operatorname{Beta}(k+1,n-k)$ respectively ($L(0)=0$, $U(n)=1$). The one-sided Clopper--Pearson bounds \cite{ClopperPearson1934} and a union bound over six endpoints with $\alpha=.05/6$, combined with \eqref{eq:finite-grid-transfer}, give simultaneous coverage at least $0.95$ for
\[
 [\max(0,L(k_-)-10^{-4}),\ \min(1,U(k_+)+10^{-4})].
\]
The decimal endpoints in the table are rounded outward; their binomial tail inequalities can be checked by integer sums with geometric remainder bounds from decreasing ratios.

\begin{table}[t]
\centering
\caption{Continuous-probability intervals for stationary OU paths at $b=0.66$ using the exact grid transitions and the deterministic time-integration bound. Each row uses $50,000$ independent grid samples. The approximation columns use a 256-node midpoint saddle, in ordinary double precision. The simultaneous coverage statement concerns the interval column.}
\label{tab:ou-continuous}
\footnotesize
\begin{tabular}{rrrrrrr}
\toprule
$d$ & $e$ & $k_-$ & $k_+$ & Interval & $P_{0,d}$ & $P_{\rm LR}$\\\midrule
32 & .0004 & 18427 & 18646 & [.36327,.37822] & .80211 & .37434\\
128 & .0002 & 9493 & 9627 & [.18557,.19690] & .28458 & .19002\\
512 & .0001 & 1428 & 1460 & [.02670,.03116] & .03607 & .03048\\\bottomrule
\end{tabular}
\end{table}

All three intervals, calculated independently of the saddle, exclude the leading approximation and contain the Lugannani--Rice values.

\medskip\noindent\textit{Rare events by tilted radial sampling.}
To estimate rarer continuous probabilities, use the trapezoidal-grid tilt. For one type set $a=e^{-\kappa h}$, $v=1-a^2$, and let $a_j,g_j$ be the weights and potential at nodes $j=0,\ldots,N$. The backward transform including node $j$ is $\exp(A_j-r_jx^2/2)$; integrating against $N(ax,v)$ gives
\[
 \begin{gathered}
 r_N=2p a_Ng_N,\qquad A_N=0,\\
 r_j=2p a_jg_j+\frac{a^2 r_{j+1}}{1+v r_{j+1}},\qquad
 A_j=A_{j+1}-\tfrac12\log(1+v r_{j+1}).
 \end{gathered}
\]
Thus the tilted initial variance is $(1+r_0)^{-1}$, and its transition from $x$ is Gaussian with mean $\widetilde a_j x$ and variance $\widetilde v_j$, where
\[
 \widetilde a_j=\frac{a}{1+v r_{j+1}},\qquad
 \widetilde v_j=\frac{v}{1+v r_{j+1}}.
\]
The variance recursion is $c_{j+1}=\widetilde a_j^2c_j+\widetilde v_j$, with $c_0=(1+r_0)^{-1}$ and $Z_N=e^{A_0}/\sqrt{1+r_0}$. The covariance-diagonal iteration is linear in the node count. The next sum of $n$ squares, conditional on $R$, has law $\widetilde v_j\chi_n^2(\widetilde a_j^2R/\widetilde v_j)$, so sampling requires only one radial coordinate per type and the likelihood is the corresponding quadrature transform.

For $d=2048,8192$ at $b=0.66$ take respectively 2048 and 4096 steps, shifts $10^{-4}$ and $2.5\times10^{-5}$, and 50,000 independent tilted draws at the grid saddle for $b$. Here $\delta_h(e)<1.1\times10^{-20}$ in both rows, and $\theta\simeq0.219$ lies in the second-moment range of Corollary~\ref{cor:is-moment}; shifted-event weights also have finite second moment because the normal coordinate is at most $R_d+de$. Table~\ref{tab:deep-ou} estimates the three grid CDFs with common samples; the exact shifted CDFs bracket $P_d(b)$ up to $\delta_h(e)$.

\begin{table}[b]
\centering
\caption{Tilted radial estimates for rare continuous OU shortfalls at $b=0.66$. Values and sampling standard errors in parentheses are in units of $10^{-5}$ for $d=2048$ and $10^{-14}$ for $d=8192$. The two right columns use the continuous Riccati saddle. The shifted CDFs determine the integration-error bracket; their Monte Carlo estimates are not exact confidence limits.}
\label{tab:deep-ou}
\footnotesize
\setlength{\tabcolsep}{3pt}
\begin{tabular}{rrrrrr}
\toprule
$d$ & $F_h(b-e)$ & $F_h(b)$ & $F_h(b+e)$ & $P_{0,d}$ & $P_{\rm LR}$\\
\midrule
2048 & 6.7796 (.0614) & 7.0851 (.0640) & 7.4047 (.0668) & 7.4430 & 7.0451\\
8192 & 1.0045 (.0132) & 1.0462 (.0137) & 1.0915 (.0143) & 1.0803 & 1.0640\\
\bottomrule
\end{tabular}
\end{table}

The larger-count probability is of order $10^{-14}$; relative second moments $5.09,9.63$ at the two counts agree in scale with $\sqrt d$ in \eqref{eq:is-relative-second-moment}. The integration shifts and sampling errors do not resolve the $d^{-1}$ relative-error order; these are examples of rare continuous probability estimation.

\medskip\noindent\textit{Reliability comparison.}
At $d=512$ compare rates $(1,4)$ and $(2,8)$ with the same powers and proportions. Proposition~\ref{prop:edgeworth} gives $y_*\simeq-0.25986$ in \eqref{eq:crossing-displacement}, predicting a crossing near $\log2+y_*/d\simeq0.69264$. Use thresholds $b_\pm=b_0+(y_*\pm1.5)/d$, 50,000 prior draws per model, 4096 steps and $e=8\times10^{-5}$. The interval procedure above with $\alpha=.05/12$ for twelve endpoints and integration allowance $10^{-4}$ gives Table~\ref{tab:ou-crossing}.

\begin{table}[t]
\centering
\caption{Continuous service CDF intervals with the same instantaneous power distributions and different temporal covariances. All six intervals have simultaneous coverage at least $0.95$ under exact grid sampling. The thresholds $b_\pm$ are rounded for display.}
\label{tab:ou-crossing}
\small
\begin{tabular}{lrr}
\toprule
Threshold & Rates $(1,4)$ & Rates $(2,8)$\\
\midrule
$0.66000$ & $[.02894,.03374]$ & $[.00560,.00794]$\\
$b_-=0.689710$ & $[.42884,.44456]$ & $[.41136,.42805]$\\
$b_+=0.695569$ & $[.55983,.57539]$ & $[.58273,.59936]$\\
\bottomrule
\end{tabular}
\end{table}

At threshold $0.66$ and target probability $0.02$, only the faster model satisfies the target. At $b_-$ its upper CDF bound is below the slower model's lower bound, with the reverse ordering at $b_+$. Thus the CDFs cross in this interval with simultaneous confidence at least $0.95$. This interval does not resolve the predicted $d^{-1}$ displacement.

\subsection{Checking the continuous Riccati calculation}

For the same two OU types (rates $(1,4)$) at $b=0.66$, independently of dense covariance calculations, solve the Riccati ODEs with cubic-spline potential on 129 equispaced points including endpoints and DOP853 integration (nodal log-residual below $2\times10^{-11}$). Refinements to 257 and 513 points leave all displayed Riccati digits unchanged. Table~\ref{tab:riccati-check} illustrates the second-order convergence of Proposition~\ref{prop:midpoint-error}.

\begin{table}[t]
\centering
\caption{Two continuous-saddle evaluation routes at $b=0.66$ for the OU model. The Riccati row uses scalar and matrix ODE integration; the other rows use $N$ midpoint nodes.}
\label{tab:riccati-check}
\small
\begin{tabular}{lrrrr}
\toprule
Method & $\theta$ & $I$ & $D$ & $V$\\\midrule
Midpoint 128 & .21888416 & .00357310532 & 1.06170335 & .14466210\\
Midpoint 256 & .21893783 & .00357395647 & 1.06171895 & .14462359\\
Midpoint 512 & .21895125 & .00357416934 & 1.06172285 & .14461396\\
Riccati & .21895573 & .00357424030 & 1.06172416 & .14461075\\\bottomrule
\end{tabular}
\end{table}

\begin{figure}[b]
\centering
\includegraphics[width=.8\linewidth]{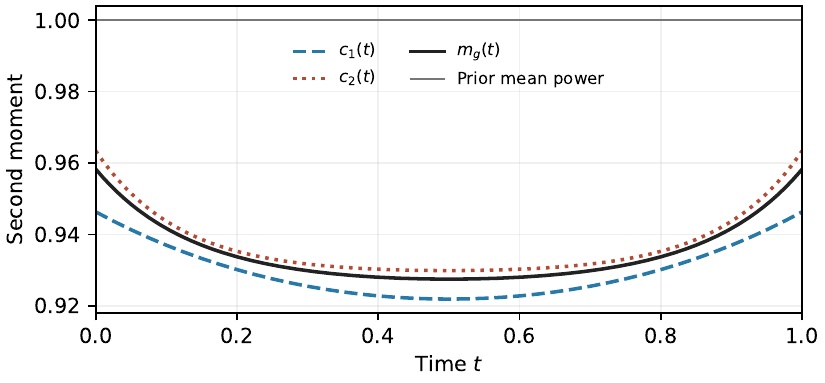}
\caption{Conditional limiting variance profiles $c_k(t)$ and the weighted aggregate $m_g(t)=(0.6c_1(t)+1.4c_2(t))/2$ for the stationary OU example at $b=0.66$. The prior mean power is one.}
\label{fig:optimal-ou-profile}
\end{figure}

The minimizing power is lowest in the interior (Figure~\ref{fig:optimal-ou-profile}): aggregate values are $0.95828$ at the endpoints and $0.92744$ at the midpoint. The stationary radial action \eqref{eq:stationary-radial-action} on the ODE solution is $0.003574240305$, agreeing with the entropy in Table~\ref{tab:riccati-check}.

For $d=2^{20}$, the Riccati value $\log P_{0,d}(0.66)\simeq-3753.16737$ is about $0.298$ below the 256-node midpoint value (about $35\%$ on the probability scale). At $r=I_H(0.66)$ the correction in \eqref{eq:quantile-expansion} adds about $2.31\times10^{-5}$ to the rate-only threshold. These are floating-point evaluations of asymptotic formulas, not bounds for actual probabilities or quantiles.

\begin{acks}[Acknowledgments]
The authors contributed equally to this work. Haoran Xu is the corresponding author.
\end{acks}

\begin{acks}[Disclosure of AI assistance]
During manuscript preparation and revision, ChatGPT and Codex assisted with language editing, intermediate mathematical derivations, and the organization and implementation of numerical calculations. The authors are responsible for the content of the manuscript.
\end{acks}

\end{document}